\documentclass{article}
\usepackage{graphicx}
\usepackage[latin2]{inputenc}
\usepackage{amsmath,amssymb,amsbsy,amsthm,amsfonts}
\usepackage{epsfig} % postscript figures
\begin{document}
%\numberwithin{equation}{section} \marginparwidth=2cm
\def\note#1{\marginpar{\small #1}}

\def\tens#1{\pmb{\mathsf{#1}}}
\def\vec#1{\boldsymbol{#1}}

\def\norm#1{\left|\!\left| #1 \right|\!\right|}
\def\fnorm#1{|\!| #1 |\!|}
\def\abs#1{\left| #1 \right|}
\def\ti{\text{I}}
\def\tii{\text{I\!I}}
\def\tiii{\text{I\!I\!I}}

\def\diver{\mathop{\mathrm{div}}\nolimits}
\def\grad{\mathop{\mathrm{grad}}\nolimits}
\def\Div{\mathop{\mathrm{Div}}\nolimits}
\def\Grad{\mathop{\mathrm{Grad}}\nolimits}

\def\tr{\mathop{\mathrm{tr}}\nolimits}
\def\cof{\mathop{\mathrm{cof}}\nolimits}
\def\det{\mathop{\mathrm{det}}\nolimits}

\def\lin{\mathop{\mathrm{span}}\nolimits}
\def\pr{\noindent \textbf{Proof: }}
\def\pp#1#2{\frac{\partial #1}{\partial #2}}
\def\dd#1#2{\frac{\d #1}{\d #2}}

\def\T{\mathcal{T}}
\def\R{\mathbb{R}}
\def\bx{\vec{x}}
\def\be{\vec{e}}
\def\bef{\vec{f}}
\def\bec{\vec{c}}
\def\bs{\vec{s}}
\def\ba{\vec{a}}
\def\bn{\vec{n}}
\def\bphi{\vec{\varphi}}
\def\btau{\vec{\tau}}
\def\bc{\vec{c}}
\def\bg{\vec{g}}

\def\bW{\vec{W}}
\def\bT{\tens{T}}
\def\bD{\tens{D}}
\def\bF{\tens{F}}
\def\bB{\tens{B}}
\def\bV{\tens{V}}
\def\bS{\tens{S}}
\def\bI{\tens{I}}
\def\bi{\vec{i}}
\def\bv{\vec{v}}
\def\bfi{\vec{\varphi}}
\def\bk{\vec{k}}
\def\b0{\vec{0}}
\def\bom{\vec{\omega}}
\def\bw{\vec{w}}
\def\p{\pi}
\def\bu{\vec{u}}

\def\ID{\mathcal{I}_{\bD}}
\def\IP{\mathcal{I}_{p}}
\def\Pn{(\mathcal{P})}
\def\Pe{(\mathcal{P}^{\eta})}
\def\Pee{(\mathcal{P}^{\varepsilon, \eta})}

\def\Ln#1{L^{#1}_{\bn}}

\def\Wn#1{W^{1,#1}_{\bn}}

\def\Lnd#1{L^{#1}_{\bn, \diver}}

\def\Wnd#1{W^{1,#1}_{\bn, \diver}}

\def\Wndm#1{W^{-1,#1}_{\bn, \diver}}

\def\Wnm#1{W^{-1,#1}_{\bn}}

\def\Lb#1{L^{#1}(\partial \Omega)}

\def\Lnt#1{L^{#1}_{\bn, \btau}}

\def\Wnt#1{W^{1,#1}_{\bn, \btau}}

\def\Lnd#1{L^{#1}_{\bn, \btau, \diver}}

\def\Wntd#1{W^{1,#1}_{\bn, \btau, \diver}}

\def\Wntdm#1{W^{-1,#1}_{\bn,\btau, \diver}}

\def\Wntm#1{W^{-1,#1}_{\bn, \btau}}

%------------------------------------------------
\newtheorem{Theorem}{Theorem}[section]
%{\theorembodyfont{\rmfamily} \newtheorem{Example}{Example}}
\newtheorem{Example}{Example}[section]
\newtheorem{Lem}{Lemma}[section]
\newtheorem{Rem}{Remark}[section]
\newtheorem{Def}{Definition}[section]
\newtheorem{Col}{Corollary}[section]
\newtheorem{Proposition}{Proposition}[section]

%--------------------------------------------------------------------------
\newcommand{\Om}{\Omega}
\newcommand{ \vit}{\hbox{\bf u}}
\newcommand{ \Vit}{\hbox{\bf U}}
\newcommand{ \vitm}{\hbox{\bf w}}
\newcommand{ \ra}{\hbox{\bf r}}
\newcommand{ \vittest }{\hbox{\bf v}}
\newcommand{ \wit}{\hbox{\bf w}}
\newcommand{ \fin}{\hfill $\square$}

\newcommand{\ZZ}{\mathbb{Z}}
\newcommand{\CC}{\mathbb{C}}
\newcommand{\NN}{\mathbb{N}}
\newcommand{\V}{\zeta}
\newcommand{\RR}{\mathbb{R}}
\newcommand{\EE}{\varepsilon}
\newcommand{\Lip}{\textnormal{Lip}}
\newcommand{\XX}{X_{t,|\textnormal{D}|}}
\newcommand{\PP}{\mathfrak{p}}
\newcommand{\VV}{\bar{v}_{\nu}}
\newcommand{\QQ}{\mathbb{Q}}
\newcommand{\HH}{\ell}
\newcommand{\MM}{\mathfrak{m}}
\newcommand{\rr}{\mathcal{R}}
\newcommand{\tore}{\mathbb{T}_3}
\newcommand{\Z}{\mathbb{Z}}
\newcommand{\N}{\mathbb{N}}

\newcommand{\F}{\overline{\boldsymbol{\tau} }}

\newcommand{\moy} {\overline {\vit} }
\newcommand{\moys} {\overline {u} }
\newcommand{\mmoy} {\overline {\wit} }
\newcommand{\g} {\nabla }
\newcommand{\G} {\Gamma }
\newcommand{\x} {{\bf x}}
\newcommand{\E} {\varepsilon}
\newcommand{\BEQ} {\begin{equation} }
\newcommand{\EEQ} {\end{equation} }
\makeatletter
\@addtoreset{equation}{section}
\renewcommand{\theequation}{\arabic{section}.\arabic{equation}}

\newcommand{\hs}{{\rm I} \! {\rm H}_s}
\newcommand{\esp} [1] { {\bf I} \! {\bf H}_{#1} }

\newcommand{\vect}[1] { \overrightarrow #1}

\newcommand{\hsd}{{\rm I} \! {\rm H}_{s+2}}

\newcommand{\HS}{{\bf I} \! {\bf H}_s}
\newcommand{\HSD}{{\bf I} \! {\bf H}_{s+2}}

\newcommand{\hh}{{\rm I} \! {\rm H}}
\newcommand{\lp}{{\rm I} \! {\rm L}_p}
\newcommand{\leb}{{\rm I} \! {\rm L}}
\newcommand{\lprime}{{\rm I} \! {\rm L}_{p'}}
\newcommand{\ldeux}{{\rm I} \! {\rm L}_2}
\newcommand{\lun}{{\rm I} \! {\rm L}_1}
\newcommand{\linf}{{\rm I} \! {\rm L}_\infty}
\newcommand{\expk}{e^{ {\rm i} \, {\bf k} \cdot \x}}
\newcommand{\proj}{{\rm I}Ę\! {\rm P}}

\renewcommand{\theenumi}{\Roman{section}.\arabic{enumi}}

\newcounter{taskcounter}[section]

\newcommand{\bib}[1]{\refstepcounter{taskcounter} {\begin{tabular}{ l p{13,5cm}} \hskip -0,2cm [\Roman{section}.\roman{taskcounter}] & {#1}
\end{tabular}}}

\renewcommand{\thetaskcounter}{\Roman{section}.\roman{taskcounter}}

\newcounter{technique}[section]

\renewcommand{\thetechnique}{\roman{section}.\roman{technique}}

\newcommand{\tech}[1]{\refstepcounter{technique} {({\roman{section}.\roman {technique}}) {\rm  #1}}}

\newcommand{\B}{\mathcal{B}}
%---------------------------------------------------------------------------

\newcommand{\diameter}{\operatorname{diameter}}

%\tableofcontents

%
%
%
%\title[  Critical and subcritical $\alpha$ models of turbulence]{Mathematical results for some $\displaystyle{\alpha}$ models of turbulence with critical and subcritical regularizations }
%
%\author[H. Ali]{Hani Ali}
%\address{IRMAR , UMR CNRS 6625, Universit\'{e} Rennes1, Campus Beaulieu, 35042 Rennes cedex, France}
%\email{hani.ali@univ-rennes1.fr}

%\author[P. Kaplicky]{Petr Kaplicky}
%\address{Charles University, Faculty of Mathematics and
%Physics, Mathematical Institute\\ Sokolovsk\'{a}~83,
%186~75~Prague~8, Czech~Republic}
%\email{Kaplicky@karlin.mff.cuni.cz}

%
%
%\keywords{turbulence model, existence, weak solution}
%\subjclass[2000]{35Q30,35Q35,76F60}

\begin{titlepage}
\title{\vskip-2cm\textbf{ Large Eddy Simulation for Turbulent  Flows with Critical Regularization}\vskip0.5cm}

\author{\begin{Large}$\hbox{Hani Ali}\thanks{IRMAR, UMR 6625,
Universit\'e Rennes 1,
Campus Beaulieu,
35042 Rennes cedex
FRANCE;
hani.ali@univ-rennes1.fr}$ \end{Large}}
\end{titlepage}
\date{}
\maketitle

\begin{abstract}
In this paper, we establish the existence of  a unique  ``regular'' weak   solution  to  the  Large Eddy Simulation (LES) models of turbulence with critical regularization. 
We first consider the critical LES for the Navier-Stokes equations and we show that its
solution converges to a solution of the Navier-Stokes  equations as the averaging radii converge
to zero. Then we extend the study to the critical LES for Magnetohydrodynamics equations.\\ 
\end{abstract}
%\maketitle
%\maketitle
%\keywords{turbulence model, existence, weak solution}
%\subjclass[2000]{35Q30,35Q35,76F60}

\medskip


%\maketitle

\section{Introduction}
Let us consider the Navier-Stokes equations in a  three dimensional torus $\mathbb{T}_3$,
\begin{align}
\diver \bv &=0, \label{nsBM}\\
\bv_{,t} + \diver ({\bv \otimes \bv}) -  \nu \Delta \bv  + \nabla p &=  {\bef},\label{nsBLM}
%\diver \overline{\bw} &=0,\\
%\alpha^{2\theta}( -\Delta)^{\theta_1} \widetilde{\bw} +   \widetilde{\bw}
%&=\bw, \quad \diver \widetilde{\bw} =0,\label{TKEbis}\\
%\alpha^{2\theta}( -\Delta)^{\theta} \overline{\varphi} +   \overline{\varphi}
%&=\varphi, \quad  \hbox{for }.\label{TKE}
%-\alpha^2 \Delta \overline{\bw} +   \overline{\bw}
%&=\bw,
\end{align}
subject to $\bv(\bx, 0) = \bv_{0}(\bx)$ and periodic boundary conditions. Here, $\bv$ is the fluid velocity field, $p$ is  the pressure, $
\bef$ is the  external body forces, $\nu$ stands for the
viscosity .\\
Equations (\ref{nsBM})-(\ref{nsBLM})
are known to be the idealized physical model to compute Newtonian fluid flows. They are
also known to be unstable in numerical simulations when the Reynolds number is high,
thus when the flow is turbulent. Therefore, numerical turbulent models are needed for real
simulations of turbulent flows. 
In many practical applications, knowing the mean characteristics of the flow by averaging techniques
is sufficient. However, averaging the nonlinear term in NSE leads to the well-known
closure problem. To be more precise, if  denotes the filtered/averaged velocity field
then the Reynolds averaged NSE (RANS)
\begin{align}
%\diver \bv &=0, \label{BM}\\
\overline{\bv}_{,t} + \diver(\overline{\bv} \otimes \overline{\bv}) -  \nu \Delta \overline{\bv} +\nabla \overline{p} + \diver \mathcal{R}(\bv,\bv) = \overline{{\bef}},\label{NOTCLOSEDBLM}
%\diver \overline{\bw} &=0,\\
%\alpha^{2\theta}( -\Delta)^{\theta_1} \widetilde{\bw} +   \widetilde{\bw}
%&=\bw, \quad \diver \widetilde{\bw} =0,\label{TKEbis}\\
%\alpha^{2\theta}( -\Delta)^{\theta} \overline{\varphi} +   \overline{\varphi}
%&=\varphi, \quad  \hbox{for }.\label{TKE}
%-\alpha^2 \Delta \overline{\bw} +   \overline{\bw}
%&=\bw,
\end{align}
where  $\mathcal{R}(\bv,\bv)=\overline{\bv \otimes \bv}-\overline{\bv} \otimes \overline{\bv} $ is the Reynolds stress tensor,
is not closed because
we cannot write it in terms of $\overline{\bv}$ alone.  The main essence of turbulence modeling is
to derive simplified, reliable and computationally realizable closure models.
%Notice that the Layton Lewandowski  model  \ref{} differs from the one introduced by Bardina et al2 where
%the following approximation is used: 
%In 1980, Bardina et al. \cite{bardina80} introduced  a particular closure model by approximating
%the Reynolds stress tensor by
%\begin{align}
%\mathcal{R}(\bv,\bv)= \overline{\overline{\bv} \otimes \overline{\bv}}-\overline{\overline{\bv}} \otimes \overline{\overline{\bv}}.
%\end{align}
In \cite{LL03} and \cite{LL06b} Layton and Lewandowski suggested an approximation of the Reynolds
stress tensor, given by
 \begin{align}
\mathcal{R}(\bv,\bv)= \overline{\overline{\bv} \otimes \overline{\bv}}-{\overline{\bv}} \otimes {\overline{\bv}}.
\end{align}
This is equivalent form to the approximation
 \begin{align}
\diver ( \overline{{\bv} \otimes {\bv}})\approx \diver (\overline{\overline{\bv} \otimes {\overline{\bv}}}).
\end{align}
Hence, Layton and Lewandowski studied the following Large Scale Model considered
as a Large Eddy Simulation (LES) model: 
\begin{align}
\diver \bw &=0, \label{llBM}\\
\bw_{,t} + \diver (\overline{\bw \otimes \bw}) -  \nu \Delta \bw + \nabla q &=  \overline{\bef},\label{llBLM}
%\diver \overline{\bw} &=0,\\
%\alpha^{2\theta}( -\Delta)^{\theta_1} \widetilde{\bw} +   \widetilde{\bw}
%&=\bw, \quad \diver \widetilde{\bw} =0,\label{TKEbis}\\
%\alpha^{2\theta}( -\Delta)^{\theta} \overline{\varphi} +   \overline{\varphi}
%&=\varphi, \quad  \hbox{for }.\label{TKE}
%-\alpha^2 \Delta \overline{\bw} +   \overline{\bw}
%&=\bw,
\end{align}
considered in $(0,T)\times \mathbb{T}_3$ and subject to $\bw(\bx, 0) = \bw_{0}(\bx)=\overline{\bv_{0}}$ and periodic boundary conditions with mean value equal to zero.
Where they denoted  $(\bw,q)$  the approximation of  $ (\overline{\bv}, \overline{p})$.\\
 The averaging operator chosen in (\ref{llBLM})  is a differential filter, \cite{germano}, \cite{GH05}, \cite{CFHOTW99b}, \cite{LL03}, \cite{dunca06}, \cite{CLT06}, that commutes with differentiation under periodic boundary conditions and is defined as follows. Let
$\alpha > 0$, given a periodic function $ \varphi \in L^2(\mathbb{T}_3)$,
define its average $  \overline{\varphi}  $  to be the unique  solution of
\begin{align}
-\alpha^{2}\Delta^{} \overline{\varphi} +   \overline{\varphi}
&=\varphi,\label{FILTRE}
\end{align}
with periodic conditions, and fields with mean value equal to zero. 

The main goal in using such a model is to filter eddies of scale less than the numerical grid
size $ \alpha$  in numerical simulations.  For a general
overview of LES models, the readers are refered  to  Berselli et al. \cite{berselli} and references cited
therein.
Notice that the Layton-Lewandowski model (\ref{llBM})-(\ref{llBLM}) differs from the one introduced by Bardina et al. \cite{bardina80} where
the following approximation of the Reynolds stress tensor is used: 
%In 1980, Bardina et al. \cite{bardina80} introduced  a particular closure model by approximating
%the Reynolds stress tensor by
\begin{align}
\mathcal{R}(\bv,\bv)= \overline{\overline{\bv} \otimes \overline{\bv}}-\overline{\overline{\bv}} \otimes \overline{\overline{\bv}}.
\end{align}

In \cite{LL03} and \cite{LL06b} Layton and Lewandowski have proved that (\ref{llBM})-(\ref{llBLM}) have a unique regular
solution . They  have aslo shown that there exists a sequence $\alpha_j $ which converges
to zero and such that the sequence $ ( \bw_{\alpha_{j}},q_{\alpha_{j}})$  converges to a distributional
solution  $(\bv,p)$ of the Navier-Stokes equations .\\
We remark that many of these  results  established in the above
cited papers have been extended  to the following three dimensional magnetohydrodynamic
equations (MHD):
\begin{align}
 \partial_t\bv-\nu_1 \Delta \bv +\diver({\bv\otimes\bv}^{})-\diver({\B\otimes\B}^{})
+\nabla p&=0,
\label{mhdIMHDV1}\\
\partial_t \B -\nu_2 \Delta \B +\diver({\bv\otimes\B}^{})-\diver({\B\otimes\bv}^{})&=0,
\label{mhdIMHDV2}\\
\displaystyle\int_{\tore} \B\,d\bx =  \int_{\tore} \bv\,d\bx =0, \quad \diver \B= \diver \B &= 0,\label{mhdIMHDV3}\\
 \B(0)=\B_{0},\; \bv(0)&=\bv_{0},\label{mhdIMHDV4}
\end{align}
 here,   $\bv$ is the fluid velocity field,  $p $ is the fluid pressure,  $\B$ is the   magnetic field, 
and $\bv_{0}$ and  $ \B_{0}$ are the corresponding initial data. 
The interested readers are referred to \cite{latr2010, LT2007}  and references cited
therein.

This paper has two main correlated objects. The first one is to study the Large Eddy Simulation  for the Navier-Stokes equations (LES for NSE)  with a general filter ${-}^{\theta}$:
\begin{align}
\diver \bw &=0, \label{BM}\\
\bw_{,t} + \diver (\overline{\bw \otimes \bw}^{\theta}) -  \nu \Delta \bw  + \nabla q&= \overline{\bef},\label{generalBLM}\\
%\diver \overline{\bw} &=0,\\
%\alpha^{2\theta}( -\Delta)^{\theta_1} \widetilde{\bw} +   \widetilde{\bw}
%&=\bw, \quad \diver \widetilde{\bw} =0,\label{TKEbis}\\
\alpha^{2\theta}( -\Delta)^{\theta} \overline{\varphi}^{\theta} +   \overline{\varphi}^{\theta}
&=\varphi,\label{vnzTKE}
%\bw_{0}&=
\end{align}
where the nonlocal operator $(-\Delta)^{\theta}$  is
defined through the Fourier transform
\begin{equation}
\widehat{(-\Delta)^{\theta}{\varphi(\bk)}}=|{\bk}|^{2\theta}\widehat{\varphi}({\bk}).
\end{equation}
%Fractionnal order Laplace operator have been employed to  another $ \alpha$ model of turbulence in  \cite{OT2007}.
Our task is to show that for $\theta \ge \frac{1}{6}$ (see Theorem \ref{TH1}), we get global in time existence of a unique weak solution $(\bw,q) $ to eqs. \eqref{BM}--\eqref{vnzTKE} such that:\\
 $(\bw,q)$ are spatially periodic with period $L$,
\begin{equation}
\label{bc1}
\int_{\tore}\bw(t,\bx) d\bx=0 \quad \textrm{ and } \int_{\tore}q(t,\bx) d\bx=0 \quad \textrm{ for } t \in [0,T),
\end{equation}
and  
\begin{equation}
\bw(0,x)=\bw_0(x)=\overline{\bv_0}^{\theta} \quad \textrm{ in } \mathbb{T}_3.
\label{ID}
\end{equation}
%Concerning  the regularized velocities $\widetilde{\bw}$, $\overline{\bw}$ we will assume that they verify the same boundary conditions as $\bw$:
% \begin{align}
%\widetilde{\bw}(t, \bx + L \vec{e}_{\vec{j}})=\widetilde{\bw}(t,\bx) \quad \textrm{ and  } \int_{\tore}\widetilde{\bw}(t,\bx) d\bx=0 \quad &&\textrm{ on } (0,T)\times \mathbb{T}_3 .\label{bctilde123}\\
%\overline{\bw}(t, \bx + L \vec{e}_{\vec{j}})=\overline{\bw}(t,\bx) \quad \textrm{ and  } \int_{\tore}\overline{\bw}(t,\bx) d\bx=0 \quad &&\textrm{ on } (0,T)\times \mathbb{T}_3,\label{bc123}
%\end{align}
We note that the value $\theta=\frac{1}{6}$ is optimal and of course the general $\alpha$ family considered here  recover the  case $\theta=1$ studied in \cite{LL06b}. The LES for NSE   with $\theta=1 $ can be also addressed as the
zeroth order Approximate Deconvolution Model referring to the family
of models in \cite{adamsstolz} . The value $ \theta=\frac{1}{6}$  is consistent with the critical regularization value needed to get existence and uniqueness to the simplified Bardina model studied in \cite{HA11}.  We note also that fractionnal order Laplace operator has been used in  another $ \alpha$ models of turbulence in  \cite{OT2007,A01,HLT10}.

The second 
object of this paper is  to study the (LES) model for Magnetohydrohynamics equation (LES for MHD)  with a general filter ${-}^{\theta}$.  Hence, we consider the following LES for MHD problem
\begin{align}
 \partial_t\bw-\nu_1 \Delta \bw +\diver(\overline{\bw\otimes\bw}^{\theta})-\diver(\overline{\vec{W}\otimes\vec{W}}^{\theta})
+\nabla q&=0,
\label{IMHDV1}\\
\partial_t \vec{W} -\nu_2 \Delta \vec{W} +\diver(\overline{\bw\otimes\vec{W}}^{\theta})-\diver(\overline{\vec{W}\otimes\bw}^{\theta})&=0,
\label{IMHDV2}\\
\displaystyle\int_{\tore} \vec{W}\,d\bx =  \int_{\tore} \bw\,d\bx =0, \quad \diver \vec{w}= \diver \vec{W} &= 0,\label{IMHDV3}\\
 \vec{W}(0)=\vec{W}_{0},\; \bw(0)&=\bw_{0},\label{IMHDV4}
\end{align}
where the boundary conditions are taken to be periodic,   and we take  as before the same filter ${-}^{\theta}$ and  $(\bw, \bW, q)$ is the approximation of  $ (\overline{\bv}, \overline{\B}, \overline{p})$ solution of the MHD Equations.\\
% Here, the unknowns are the averaging fluid velocity field $\bw(t,\vec{x})$, the averaging fluid pressure $q(t,\vec{x})$, and the averaging  magnetic field $\vec{W}(t,\vec{x})$.  Note that when $\alpha=0$, we formally retrieve the  MHD equations . 
The case when $\theta=1$ is studied in \cite{latr2010} where the authors gived a  mathematical description of the problem, 
performed the numerical analysis of the model and verified their physical
fidelity.\\
In this paper, we show   that for $\theta \ge \frac{1}{6}$ (see Theorem \ref{existensssIMHDV}), we get global in time existence of a unique weak solution $(\bw, \bW, q) $ to eqs. \eqref{IMHDV1}--\eqref{IMHDV4}.
Let us mention that the idea to consider the LES for MHD with critical regularization  is a new feature for the present work.
The Approximate Deconvolution Model for  Navier-Stokes equations with $ \theta > \frac{3}{4}$ is studied in \cite{bresslilewandowski}  and the  Approximate Deconvolution Model for Magnetohydrodynamics equations with $\theta=1$ is studied in \cite{latr210}. 
 As mentioned  in \cite{bresslilewandowski}   the value $ \theta > \frac{3}{4}$  is not optimal in order to prove the existence and the uniqueness of the solutions in the deconvolution case.\\ 
 \textit{ ``The exponent "3/4" looks like a "critical exponent". We conjecture that we
can get an existence and uniqueness result for lower exponents, but concerning the convergence towards the mean Navier-Stokes equations, we think that it is the best exponent, but
this question remains an open one."}\\
Notice however that  unlike the LES case the value $\theta=\frac{1}{6}$ is not sufficent to get the uniqueness in the deconvolution case. 
Based on this work, we will study in a forthcoming paper the Approximate Deconvolution Model for both Navier-Stokes equations and  Magnetohydrodynamics equations with critical regularizations. \\
Finally, one may ask questions about the relation between the regularization parameter $\theta$ and  the model consistency errors. These  questions are adressed in \cite{LL06} where $\theta=1$. 
Therefore, the issue is to find the relation between the model consistency errors and the regularization parameter $\theta$.\\

This  paper is organized as follows. In section 2 we prove the global existence and uniqueness of the solution for the LES for NSE with critical regularization.  We also prove that the solution $({\bw},q)$ of the LES for NSE   converges in
some sense to a solution of the Navier-Stokes equations  when $\alpha$  goes to zero. Section 3 treats the
questions of global existence, uniqueness and convegence for the LES for MHD with critical regularization.

\section{ The Critical LES for NSE }

\ Before formulating the main results of this paper, we fix notation of function spaces that we shall  employ.\\ 
We denote by $L^p(\tore)$ and $W^{s,p}(\tore)$, $s \ge -1 $,  $1 \le p \le \infty$ the usual Lebesgue and Sobolev spaces over $\tore$, and the Bochner spaces $C(0,T;X), L^p(0,T;X)$ are defined in the standard way. 
%For the velocity field, we define
In addition we introduce 
\begin{align*}
\dot{W}^{s,p}_{\diver}&=\left\{ \bw \in W^{s,p}(\tore)^3; \; \int_{\tore}\bw=0 ; \; \diver \bw
=0 \textrm{ in }  \tore \right\}.
%W^{-1,p'}_{}&=\left(W^{1,p}_{} \right)^{'}, \quad
%W^{-1,p'}_{\diver}=\left(W^{1,p}_{\diver} \right)^{'},\\
%L^2_{}&= \overline{W^{1,2}_{}}^{\|\, \|_{2}}, \quad 
%L^2_{\diver}= \overline{W^{1,2}_{\diver}}^{\|\, \|_{2}}.
\end{align*}
% Let us mention that by using Poincar\'e inequality we have 
%  \begin{align}
% \|\vec{v} \|_{s,2} \approx \| \widetilde{\vec{v}}\|_{s+2\theta_1,2} \approx \| \overline{\vec{v}} \|_{s+2\theta_2,2}. 
% \end{align}
%For $s  \in \R$,
%let $\vec{I}=\left\lbrace \bk  \ \hbox{such  that}  \   \displaystyle \bk=\frac{2\pi \ba}{\vec{L}},\  \ba \in \Z^{3} , \ba \neq 0 \right\rbrace ,$\\
%the usual Sobolev spaces $W^{s,2}$  with zero space average are classicaly characterized in terms of Fourier series  
%$$W^{s,2}=\{\bw=\sum_{\bk \in \vec{I}}\hat{\bw}_{\bk}\exp{\left\lbrace   i\bk\cdot \bx\right\rbrace }, \
%\| \bw \|_{{s,2}}^{2}< \infty          \},$$
%where 
%$$\| \bw \|_{{s,2}}^{2}= \sum_{\bk \in \vec{I}}|\bk|^{2s} |\hat{\bw}_{\bk}|^{2}.$$
%Therefore,  the space $W^{s,2}_{\diver}$ is  
%$$W^{s,2}_{\diver}=\{\bw=\sum_{\bk \in \vec{I}}\hat{\bw}_{\bk}\exp{\left\lbrace   i\bk\cdot \bx\right\rbrace }, \ 
%\ \vec{k}\cdot\hat{\bw}_{\bk}=0, \
%\| \bw \|_{{s,2}}^{2}< \infty          \}.$$
We present our main results, restricting ourselves  to the  critical case $\theta = \frac{1}{6}$,
 and for simplicity we drop some indices of $\theta$ so sometimes we will write   ``${\overline{\varphi}}^{}$'' instead of ``${\overline{\varphi}}^{\theta} $'', expecting that no confusion will occur.
\subsection{Existence and uniqueness results for the LES for NSE}
\begin{Theorem}
Assume that $\theta = \frac{1}{6}$. Let  $\overline{\bef} \in L^{2}
(0,T;W^{-\frac{5}{6},2}_{})$ be a divergence free function and  $\bw_0 \in W^{\frac{1}{6},2}_{\diver}$. Then there exist $(\bw, q)$ a unique ``regular" weak  solution to \eqref{BM}--\eqref{ID}  such that
\begin{align}
\bw &\in \mathcal{C}_{}(0,T;\dot{W}^{\frac{1}{6},2}_{\diver}) \cap
L^2(0,T;\dot{W}^{1+\frac{1}{6},2}_{\diver}),\label{bv12}\\
\bw_{,t}&\in   L^{2}(0,T;W^{-\frac{5}{6},2}_{}),
\label{bvt}\\
q&\in L^{2}(0,T;W^{\frac{1}{6},2}(\tore))
\label{psp}.
\end{align}
fulfill
\begin{equation}
\begin{split}
\int_0^T \langle \bw_{,t}, \bfi \rangle  -  (\overline{\bw \otimes \bw}, \nabla
\bfi) +  \nu(
\nabla \bw, \nabla \bfi  )\; dt =\int_0^T  \langle \bef, \bfi \rangle \; dt\\
\qquad \textrm{ for all } \bfi\in L^{2}(0,T; \dot{W}^{\frac{5}{6},2}_{\diver}),
\end{split}\label{weak1}
\end{equation}
Moreover,
\begin{equation}
\bw(0)= \bw_0.
\label{intiale}
\end{equation}
 %and the initial conditions are attained in the following sense
%\begin{equation}
%\lim_{t\to 0+}\left(\|\bw(t)-\bw_0\|_2^2 \right)
%=0. \label{inca}
%\end{equation}
 \label{TH1}
\end{Theorem}
\begin{Rem}
The notion of ``regular weak solution" is introduced in \cite{bresslilewandowski}. 
Here, we use the name ``regular"  for the weak  solution since the weak solution is unique and the velocity part of the solution  $\bw$ does not develop a finite time singularity. 
%$\in  \mathcal{C}_{}(0,T;\dot{W}^{\frac{1}{6},2}_{\diver})$
\end{Rem}
\begin{Rem}
 Once existence and  uniqueness in the large of a weak solution to the model  \eqref{BM}--\eqref{ID} with critical regularization is known. Further theoretical
properties of the model with critical and subcritical regularizations can then be developed. These are
currently under study by the author and will be presented in a subsequent report.
\end{Rem}
\textbf{Proof of Theorem \ref{TH1}.}
The proof of Theorem \ref{TH1} follows the classical scheme. We start by constructing
approximated solutions $(\bv^N, p^N)$ via Galerkin method. Then we seek for a priori
estimates that are uniform with respect to $N$. Next, we passe to the limit in the
equations after having used compactness properties. Finaly we show that the solution
we constructed is unique thanks to Gronwall's lemma.\\

\textbf{Step 1}(Galerkin approximation).
Consider a sequence $\left\{ \bfi^{r} \right\}_{r=1}^{\infty}$ consisting of $L^2$-orthonormal  and $W^{1,2}$-orthogonal eigenvectors of the Stokes problem subjected to the
space periodic conditions. 
%\begin{align}
%-\Delta \bfi^{r} = |r|^2 \bfi^{r},\hbox{ in } \tore \quad \int_{\tore}\bfi^{r}=0,  \quad {\bfi^{r}}( \bx + L \vec{e}_{\vec{j}})={\bfi^{r}}(\bx) \hbox{ for all } {r} \in \Z^3\setminus\{0\} .
%\end{align}
We note that this sequence forms a hilbertian basis of $L^2$.\\
%Stokes operator subjected to the space periodic conditions (\ref{bc1}) and 
We set 
\begin{equation}
\bw^N(t,\bx)=\sum_{{r=1}}^{N}\bc_{r}^N (t) \bfi^{r}(\bx), \hbox{ and } q^N(t,\bx)=\sum_{{|\vec{k}|=1}}^{N}q_{\vec{k}}^N (t) e^{i \vec{k} \cdot \bx}.
\end{equation}
% such that $\bk \cdot  \bc_{r}^N = 0  $ for all $\bk \in \Z^3\setminus\{0\}$ and $\left(\bc_{r}^N \right)^{*}=\bc_{-r}^N  $, where  $ \left(\bc_{r}^N \right)^{*}$ denote the complex conjugate $\bc_{r}^N.  $
% Thus due of (\ref{TKEbis}) and (\ref{TKE}) we have 
%\begin{equation}
% \widetilde{\bw}^N(t,\bx)=\sum_{{|r|=1}}^{N}\widetilde{\bc}_{r}^N (t) \bfi^{r}(\bx) \textrm{ and }  \overline{\bw}^N(t,\bx)=\sum_{{|r|=1}}^{N}\overline{\bc}_{r}^N (t) \bfi^{r}(\bx), 
%\end{equation}
% where 
%\begin{equation}
% \widetilde{\bc}_{r}^N = \frac{{\bc}_{r}^N }{1+\alpha^{2\theta_1}|r|^{2\theta_1}} \textrm{ and }  \overline{\bc}_{r}^N = \frac{{\bc}_{r}^N }{1+\alpha^{2\theta_2}|r|^{2\theta_2}}, 
%\end{equation}
%for all $\bk \in \Z^3\setminus\{0\}$.\\
We look for $(\bw^N(t,\bx), q^N(t,\bx)) $ that are determined through the system of equations 
\begin{equation}
\begin{split}
 \left( \bw_{,t}^N, \bfi^{r} \right)  -  (\overline{\bw^N \otimes \bw^N}, \nabla
\bfi^{r}) +  \nu(
\nabla \bw^N, \nabla \bfi^{r}  )\; 
= \langle \bef, \bfi^{r} \rangle \; ,
\qquad {r=1},2,...,N,
\end{split}\label{weak1galerkine}
\end{equation}
and 
%\begin{equation}
%\begin{split}
% q^N = - \sum_{i,j}\partial_i\partial_j \Delta^{-1}(\Pi^N(\overline{{v}_i^N {v}_j^N)})= - \sum_{i,j} R_{ij}(\Pi^N(\widetilde{v}_i^N\overline{v}_j^N)),
%\end{split}\label{pressuregalerkine}
%\end{equation}
\begin{equation}
\begin{split}
\displaystyle \Delta^{} q^N=  -\diver \diver  \left(\Pi^N(\overline{{\bw}^N \otimes {\bw}^N}) \right).
\end{split}\label{pressuregalerkine}
\end{equation}
Where the projector $ \displaystyle \Pi^N $ assign to any Fourier series $\displaystyle \sum_{\bk \in \Z^3\setminus\{0\}} \vec{g}_{\bk} e^{i\bk \cdot\bx} $ its N-dimensional part, i.e. 
 $\displaystyle \sum_{\bk \in \Z^3\setminus\{0\}, |\bk| \le N} \vec{g}_{\bk} e^{i\bk \cdot\bx}. $
 %we deduce from the elliptic equation (\ref{pressuregalerkine}) that 
 % after dropping indices of $N$ for shortness, that 
  % \begin{equation}
%\begin{split}
% p^N = - \sum_{i,j}\partial_i\partial_j \Delta^{-1}(\Pi^N(\widetilde{v}_i^N\overline{v}_j^N))= - \sum_{i,j} R_{ij}(\Pi^N(\widetilde{v}_i^N\overline{v}_j^N)),
%\end{split}\label{pressurepseudo}
%\end{equation}
% and $R_{ij}$ is the Riez operator defined through the Fourier transform by
% \begin{equation}
%\begin{split}
%  \widehat{R_{ij}(u)}=\frac{k_i k_j}{|\bk|^2}\widehat{u({\bk})}, \quad \hbox{ for all } \bk \in \Z^3\setminus\{0\}.
%  \end{split}\label{pressurepseudo2}
%\end{equation}

% The Laplace operator $\displaystyle \Delta^{-1}$ can be  defined through the Fourier transform hence $\displaystyle (\Delta^{-1}\diver \diver)$ may be viewed as a pseudofifferentiel operator 
%  \begin{equation}
%\begin{split}
%   ((\Delta^{-1}\diver \diver){\bw})_{\bk}=\frac{\bk^2}{|\bk|^2}\bw_{\bk}.
%\end{split}\label{pressurepseudo}
%\end{equation}
Moreover we require that $\bw^N $ satisfies the following initial condition
\begin{equation}
 \label{initial Galerkine}
\bw^N(0,.)= \bw^N_0= \sum_{r=1}^{N}\bc_0^N  \bfi^{r}(\bx),
\end{equation}
and 
\begin{equation}
 \label{initial2 Galerkine}
\bw^N_0 \rightarrow \bw_0  \quad \textrm{ strongly  in }  W^{\frac{1}{6},2}(\tore)^3  \quad \textrm{ when } N \rightarrow \infty .
\end{equation}
%Where the initial condition $ \overline{\bw}^N_0$ is deduced from $\bw^N_0$ through the relation (\ref{TKE}).
%\begin{equation}
% \label{initialbar Galerkine}
%  \overline{\bw}^N_0   + \alpha^{2\theta}(-\Delta)^{\theta} \overline{\bw}^N_0=\bw^N_0
%\end{equation}

The classical Caratheodory theory \cite{Wa70}  then implies the short-time existence of solutions 
to (\ref{weak1galerkine})-(\ref{pressuregalerkine}).  Next we derive  estimate on $\bc^N$ that is uniform w.r.t. $N$.
These estimates then imply that the  solution of  (\ref{weak1galerkine})-(\ref{pressuregalerkine}) constructed on a short time interval $[0, T^N[ $ exists for all $t \in [0, T]$.\\

\textbf{Step 2} (A priori estimates)
Multilplying the $r$th equation in (\ref{weak1galerkine}) with $\alpha^{2\theta}|\bk|^{2\theta}{\bc}^N_{r}(t)+{\bc}^N_{r}(t)$, summing over ${r=1},2,...,N$, integrating over time from $0$ to $t$ and using the following identities
 
\begin{equation}
 \label{divergencfreebar1}
\left({\bw}^N_{,t},  {\bw}^N_{}+ \alpha^{\frac{1}{3}}(-\Delta)^\frac{1}{6}{\bw}^N_{} \right)=
\frac{1}{2}\frac{d}{dt}\|{\bw}^N \|_{2}^2 + \frac{1}{2}\frac{d}{dt}\|{\bw}^N \|_{\frac{1}{6},2}^2,
\end{equation}
\begin{equation}
 \label{divergencfreebar2}
\left(-\Delta{\bw}^N,  {\bw}^N_{}+ \alpha^{\frac{1}{3}}(-\Delta)^\frac{1}{6}{\bw}^N_{} \right)=
\|{\bw}^N \|_{1,2}^2 + \|{\bw}^N \|_{1+\frac{1}{6},2}^2,
\end{equation}
\begin{equation}
 \label{divergencfreebar3}
\langle \overline{\bef}, {\bw}^N_{}+ \alpha^{\frac{1}{3}}(-\Delta)^\frac{1}{6}{\bw}^N_{} \rangle= \langle \bef, {\bw}^N \rangle,
\end{equation}
and
\begin{equation}
 \label{divergencfreebar}
 \begin{array}{lll}
\left(\overline{{\bw}^N \otimes {\bw}^N}, \nabla ( {\bw}^N_{}+ \alpha^{\frac{1}{3}}(-\Delta)^\frac{1}{6}{\bw}^N_{} )\right)&=\left(  \bw^N \otimes {\bw}^N, \nabla {{\bw}^N}\right)\\
&=-\left(  \diver {\bw}^N,  \frac{|{\bw}^N|^2}{2}\right)=0
\end{array}
\end{equation}
leads to the a priori estimates
\begin{equation}
\begin{array}{llll}
 \label{apriori1}
\displaystyle \frac{1}{2}\left(\|{\bw}^N \|_{2}^2 + \|{\bw}^N \|_{\frac{1}{6},2}^2\right)+ \displaystyle \nu\int_{0}^{t}\left(\|{\bw}^N \|_{1,2}^2  + \|{\bw}^N \|_{1+\frac{1}{6},2}^2 \right) \ ds\\
 \quad \quad \quad = \displaystyle \int_{0}^{t} \langle \bef, {\bw}^N \rangle \ ds 
 +\displaystyle \frac{1}{2}\left(\|{\bw}_0 \|_{2}^2 + \|{\bw}_0 \|_{\frac{1}{6},2}^2\right).
\end{array}
\end{equation}
 Using the duality norm comined with Young inequality we conclude from eqs. (\ref{apriori1}) that 
 \begin{equation}
 \label{apriori12}
\sup_{t \in [0,T^N[}\|{\bw}^N \|_{2}^2 + \sup_{t \in [0,T^N[}\|{\bw}^N \|_{\frac{1}{6},2}^2+ \nu\int_{0}^{t}\left(\|{\bw}^N \|_{1,2}^2  + \|{\bw}^N \|_{1+\frac{1}{6},2}^2 \right) \ ds  \le C
\end{equation}
that immediately implies that the existence time is independent of $N$ and it is possible to take $T=T^N$.\\ 
We deduce from (\ref{apriori12}) that 
\begin{equation}
\label{vbar1}
 {\bw}^N \in L^{\infty}(0,T ; \dot{W}^{\frac{1}{6},2}_{\diver}) \cap L^{2}(0,T ; W^{1+\frac{1}{6},2}(\tore)^3).
 \end{equation}
  From  ({\ref{vbar1}}) and (\ref{TKE})  it follows   that  
  \begin{equation}
\label{vtilde1}
\overline{{\bw}^N \otimes {\bw}^N} \in L^{2}(0,T ; W^{\frac{1}{6},2}_{}).  
 \end{equation}
Consequently from the elliptic theory eqs (\ref{pressuregalerkine}) implies that 
 \begin{equation}
\label{vbarvbarpressure}
\int_{0}^{T}\|p^N\|_{\frac{1}{6},2}^2 dt < K. 
\end{equation}
From eqs. (\ref{weak1galerkine}), (\ref{vbar1}), (\ref{vtilde1}) and (\ref{vbarvbarpressure}) we also obtain that 
 \begin{equation}
\label{vtemps}
\int_{0}^{T}  \|\bw^N_{,t}\|_{-\frac{5}{6},2}^2  dt < K. 
\end{equation}

\textbf{Step 3} (Limit $N \rightarrow \infty$) It follows from the estimates (\ref{vbar1})-(\ref{vtemps}) and the Aubin-Lions compactness lemma
(see \cite{sim87} for example) that there are a  not relabeled  subsequence of $(\bw^N, q^N)$  and a couple $(\bw, q)$ such that

\begin{align}
\bw^N &\rightharpoonup^* \bw &&\textrm{weakly$^*$ in } L^{\infty}
(0,T;W^{\frac{1}{6},2}), \label{c122}\\
\bw^N &\rightharpoonup \bw &&\textrm{weakly in }
L^2(0,T;W^{1+\frac{1}{6},2}_{}), \label{c22}\\
\overline{{\bw}^N \otimes {\bw}^N} &\rightharpoonup \overline{\bw\otimes {\bw}} &&\textrm{weakly in }
L^2(0,T;W^{\frac{1}{6},2}_{}), \label{c22prime}\\
\bw^N_{,t}&\rightharpoonup \bw_{,t} &&\textrm{weakly in } L^{2}
(0,T;W^{-\frac{5}{6},2}_{}),
\label{c322}\\
q^N&\rightharpoonup q &&\textrm{weakly in } L^{2}(0,T;W^{\frac{1}{6},2}
(\tore)), \label{c32}\\
%\bw^N &\rightharpoonup \bw &&\textrm{weakly in } L^{\frac{8}{3}}
%(0,T;L^{\frac{8}{3}}(\partial \Omega)^3). \label{c5.12}\\
\bw^N &\rightarrow \bw &&\textrm{strongly in  }
L^2(0,T;W^{s,2}(\tore)^3), s < 1+\frac{1}{6} \label{c83ici}\\
\overline{{\bw}^N \otimes {\bw}^N} &\rightarrow \overline{{\bw} \otimes \bw} &&\textrm{strongly in  }
L^2(0,T;W^{r,2}(\tore)^3), r < \frac{1}{6}.\label{c82int}
%\overline{\bw}^N &\rightarrow \overline{\bw} &&\textrm{strongly in  }
%L^q(0,T;L^q(\tore)^3) \textrm{  for all } q< 5 .\label{c82pppp}
\end{align}
%
%By a standard interpolation argument we have  
%
%\begin{align}
%\bw^N & \in  L^{\frac{10}{3}}(0,T;L^{\frac{10}
%{3}}(\tore)^3), \label{intrpolv}\\
%\widetilde{\bw}^N & \in  L^{\frac{10}{3-4\theta_1}}(0,T;L^{\frac{10}
%{3-4\theta_1}}(\tore)^3), \label{interpvtilde}\\
%\overline{\bw}^N & \in  L^{\frac{10}{3-4\theta_2}}(0,T;L^{\frac{10}
%{3-4\theta_2}}(\tore)^3). \label{interpvbar}
%\end{align}
%Thus  from (\ref{intrpolv})-(\ref{interpvbar}) and (\ref{c83ici})-(\ref{c82int})
%we obatin
%\begin{align}
%\bw^N &\rightarrow \bw &&\textrm{strongly in  }
%L^{q_1}(0,T;L^{q_1}(\tore)^3) \textrm{ for all } q_1<\frac{10}{3},\label{convint1}\\
%\widetilde{\bw}^N &\rightarrow \widetilde{\bw} &&\textrm{strongly in  }
%L^{q_2}(0,T;L^{q_2}(\tore)^3) \textrm{ for all } q_2<\frac{10}{3-4\theta_1},\label{convint2}\\
%\overline{\bw}^N &\rightarrow \overline{\bw} &&\textrm{strongly in  }
%L^{q_3}(0,T;L^{q_3}(\tore)^3) \textrm{ for all } q_3<\frac{10}{3-4\theta_2},\label{convint3}
%\end{align}
%Since $ q_2<\frac{10}{3-4\theta_1}$, $q_3<\frac{10}{3-4\theta_1} $ and $2\theta_1 + \theta_2 = \frac{1}{2} $  the application of H$\ddot{o}$lder's inequality implies that
%\begin{align}
% \widetilde{\bw} \otimes \overline{\bw} \in 
%L^{q}(0,T;L^{q}(\tore)^{3\times 3}) \textrm{ where } q \ge {2},\label{vtimesv}
%\end{align}
%%These limits are sufficient to pass to the limit (from Galerkin approximations (17) to a weak formulation
%%of Eq. (10)) in the terms involving time derivative, the pressure and due to Eq. (27) also the convective
%%term.

The above established convergences are clearly sufficient for taking the limit in (\ref{weak1galerkine})  and for concluding that  the velocity part  $ \bw$  satisfy (\ref{weak1}). 
Moreover, 
from (\ref{c22}) and (\ref{c322})
 one we can deduce by a classical argument ( see in \cite{A01})   that 
 \begin{equation}
 \bw \in  \mathcal{C}(0,T;W^{\frac{1}{6},2}).
\end{equation}
Furthermore, from  the  strong continuty of $\bw$ with respect to the time with value in $W^{\frac{1}{6},2}$  we deduce   that $\bw(0)=\bw_0$.\\
Let us mention also that ${\bw}+ \alpha^{\frac{1}{3}}(-\Delta)^\frac{1}{6}{\bw} $ is a possible  test in the weak formlation (\ref{weak1}). Thus $ {\bw}$ verifies for all $t \in [0,T]$  the follwing equality   
\begin{equation}
\begin{array}{lll}
 \label{apriori12leray}
\displaystyle \left(\|{\bw} (t)\|_{2}^2 + \|{\bw}(t)\|_{\frac{1}{6},2}^2\right)+\displaystyle 2\nu\int_{0}^{t}\left(\|{\bw} \|_{1,2}^2  + \|{\bw} \|_{1+\frac{1}{6},2}^2 \right)ds \\
\quad \quad \quad = \displaystyle 2\int_{0}^{t}\langle \bef, {\bw} \rangle ds + \left(\|{\bw}_0 \|_{2}^2 + \|{\bw}_0\|_{\frac{1}{6},2}^2\right).
\end{array}
\end{equation}

 \textbf{Step 5} (Uniqueness)
Since the pressure part of the solution is uniquely determined by the velocity part it remain to show the uniqueness to the velocity.

Next, we will show the continuous dependence of the  solutions on the initial data and in particular the uniqueness.\\
Let $({\bw_1,p_1})$ and $({\bw_2,p_2})$   any two solutions of (\ref{BM})-(\ref{TKE}) on the interval $[0,T]$, with initial values $\bw_1(0)$ and $\bw_2(0)$. Let us denote by  $\delta \vec{w}_{} =\bw_2-\bw_1$.
 We subtract the equation for $\bw_1$ from the equation for $\bw_2$ and test it with $\delta \vec{w}$.
 We get using successively the relation (\ref{TKE}), the fact that the averaging operator commutes with differentiation under periodic boundary conditions, the norm duality, Young inequality and Sobolev embedding theorem:
 \begin{equation}
 \begin{split}
  \displaystyle 
  \|{\delta \vec{w}}_{,t}\|_{2}^{2} + \alpha^{\frac{1}{3}}\|{\delta \vec{w}}_{,t}\|_{\frac{1}{6},2}^{2} +\nu \|\nabla{\delta \vec{w}}_{}\|_{2}^{2}+ \alpha^{\frac{1}{3}}\|\nabla{\delta \vec{w}}_{}\|_{\frac{1}{6},2}^{2} \\
   \le \displaystyle (\overline{\bw_{2}\otimes{\bw}_{2}}-\overline{\bw_{1}\otimes{\bw}_{1}}, \nabla ({\delta \vec{w}} + \alpha^{\frac{1}{3}} (-\Delta)^{\frac{1}{6}}{\delta \vec{w}}_{}))
    \\   \le \displaystyle
    ( {\bw_{2}\otimes{\bw}_{2}}-{\bw_{1}\otimes{\bw}_{1}}, \nabla {\delta \vec{w}} )\\
\le \displaystyle \frac{4}{\nu} \|{\delta \vec{w}}\otimes{\bw}_{1} \|_{-\frac{1}{6},2}^{2}    \\
\quad \quad \le \displaystyle \frac{4}{\nu} \|{\delta \vec{w}}\|_{\frac{1}{6},2}^{2} \|{\bw}_{1} \|^{2}_{1+\frac{1}{6}}.
\end{split}
\end{equation}
 % \begin{equation}
% \begin{split}
%  \displaystyle 
%\|{\delta \vec{w}}_{,t}\|_{2}^{2} +\nu \|\nabla{\delta \vec{w}}_{}\|_{2}^{2} & \le \displaystyle
%\frac{4}{\nu} \|\widetilde{\delta \vec{w}}
%\overline{\bw}_{1} \|_{2}^{2} +  \frac{4}{\nu} \|\widetilde{\vec{v}_2}
%\overline{\delta \vec{w}} \|_{2}^{2}    \\
%&\le \displaystyle \frac{4}{\nu} \|\widetilde{\delta \vec{w}}\|_{\frac{1}{2}-2\theta_2,2}^{2}
%\|\overline{\bw}_{1} \|^{2}_{1+2\theta_2}+  \displaystyle\frac{4}{\nu} \|\overline{\delta \vec{w}}\|_{\frac{1}{2}-2\theta_1,2}^{2}
%\|\widetilde{\bw}_{2} \|^{2}_{1+2\theta_1}  \\ &\le \displaystyle \frac{1}{\alpha^{{2\theta_1}+{2\theta_2}}}\frac{4}{\nu} \|{\delta \vec{w}}\|_{2}^{2}
%\left(\|{\bw}_{1} \|^{2}_{1,2} + 
%\|{\bw}_{2} \|^{2}_{1,2}\right)
%\end{split}
%\end{equation}
%\begin{equation}
% \begin{split}
%  \displaystyle 
%\|{\delta \vec{w}}_{,t}\|_{2}^{2} +\nu \|\nabla{\delta \vec{w}}_{}\|_{2}^{2} & \le \displaystyle
%\frac{4}{\nu} \|\widetilde{\delta \vec{w}}
%\overline{\bw}_{1} \|_{2}^{2} +  \frac{4}{\nu} \|\widetilde{\vec{v}_2}
%\overline{\delta \vec{w}} \|_{2}^{2}    \\
%&\le \displaystyle \frac{4}{\nu} \|\widetilde{\delta \vec{w}}\|_{\frac{3}{4}+2\theta_1-\theta_2,2}^{2}
%\|\overline{\bw}_{1} \|^{2}_{\frac{3}{4}-2\theta_1+\theta_2}+  \displaystyle\frac{4}{\nu} \|\overline{\delta \vec{w}}\|_{\frac{1}{2}-2\theta_1+\theta_2,2}^{2}
%\|\widetilde{\bw}_{2} \|^{2}_{1+2\theta_1-\theta_2}  \\ &\le \displaystyle \frac{1}{\alpha^{{2\theta_1}+{2\theta_2}}}\frac{4}{\nu} \|{\delta \vec{w}}\|_{2}^{2}
%\left(\|{\bw}_{1} \|^{2}_{1,2} + 
%\|{\bw}_{2} \|^{2}_{1,2}\right)
%\end{split}
%\end{equation}
Using Gronwall's inequality we conclude 
the continuous dependence of the solutions on the inital data in the $L^{\infty}([0,T],W^{\frac{1}{6},2}_{})$  norm. In particular, if ${\delta \vec{w}}^{}_{0}=0$ then ${\delta \vec{w}}=0$ and the solutions are unique for all $t \in [0,T] .$  Since $T>0$ is arbitrary this solution may be uniquely extended for all time.\\ 
This finish the proof of Theorem \ref{TH1}.\\

\subsection{Limit consistensy for the critical LES for NSE }
In this section, we  take the limit $\alpha \rightarrow 0$ in order to show the following result:

\begin{Theorem}
\label{1deuxieme}
 Let $(\bw_{\alpha}, q_{\alpha})$ be the solution of \eqref{BM}--\eqref{vnzTKE}  for a fixed $\alpha$.
 There is a subsequence $\alpha_j$  such that
$(\bw_{\alpha_j}, q_{\alpha_j}) \rightarrow (\bv, p)$ as $j \rightarrow \infty $
where
$ (\bv, p)  \in L^{\infty}
([0,T];L^2({\tore})^3)\cap L^2([0,T];\dot{W}^{1,2}_{\diver}) \times L^{\frac{5}{3}}([0,T];L^{\frac{5}{3}}(\tore))$
is a weak solution of the Navier-Stokes equations  with periodic boundary conditions and zero mean value constraint.\\
 The sequence  $ \bw_{\alpha_j}$  converges strongly
to $\bv$ in the space $L^p([0,T];L^p(\tore)^3)$ for all  $2 \le  p<\frac{10}{3}, $  and weakly in $ L^{r}(0,T; L^{\frac{6r}{3r-4}}(\tore)^3) \textrm { for all } r \ge 2, $ while the sequence
 $  q_{\alpha_j} $  converges strongly to $p$ in the space $ L^p([0,T];L^p(\tore)) \textrm{ for all } \frac{4}{3}\le  p<\frac{5}{3}, $ and weakly  in the space  $L^{\frac{r}{2}}(0,T; L^{\frac{3r}{3r-4}}(\tore)) \textrm { for all } r \ge 2.$
\end{Theorem}
Before proving Theorem \ref{1deuxieme}, we first record the following three Lemmas.

%On commence par le Lemme suivant:
%We start by the following Lemma: 
% On prend $2\phi\vec{u}_{\alpha}$ comme une fonction test dans (\ref{alpha ns}). On note que la condition $\theta \ge 1/4 $ assure que  toutes les termes obtenues sont bien définies, en particulier l'intégrale $$\displaystyle 2\int_{0}^{T} \int_{\tore} \overline{\vec{u}_{\alpha}}\nabla \vec{u}_{\alpha} \cdot \vec{u}_{\alpha} \phi   \  d\vec{x} dt $$ est finie en utilisant le fait que $ \overline{\vec{u}_{\alpha}}\nabla \vec{u}_{\alpha} \in L^2([0,T];\vec{H}^{-1})$ et $2\phi\vec{u}_{\alpha} \in L^2([0,T];\vec{H}^{1})$.\\ 
%Maintenant, des intégrations par parties en utilisant le fait que $\phi(T,\cdot)=\phi(0,\cdot)=0$ et l'identité suivante : 
%\begin{equation}
%2\int_{\tore} \overline{\vec{u}_{\alpha}}\nabla \vec{u}_{\alpha} \cdot \vec{u}_{\alpha} \phi   \  d\vec{x} = \int_{\tore} \overline{\vec{u}_{\alpha}} |\vec{u}_{\alpha}|^2   \cdot \nabla \phi   \  d\vec{x}
%\end{equation}
%montrent que  $(\vec{u}_{\alpha}, p_{\alpha})$ vérifie bien (\ref{local alpha}).\\
%
%Pour prendre la limite $\alpha \rightarrow 0$ on a besoin de montrer que pour tout $\bu_\alpha \in L^p(0,T;L^p(\tore)^3) $, $2 \le p<10/3.$  on a :  
%\begin{align}
%\overline{\bu_\alpha}  &\rightarrow \bu &&\textrm{strongly in  }
%L^p(0,T;L^p(\tore)^3) \textrm{ pour tout }  2 \le p<10/3.
%\end{align}
%Cela est le but de ces deux Lemmes suivants. 
\begin{Lem}
\label{fourierdiscret}
Let $ \theta \in\RR^+,  0 \le \beta \le 2\theta , s\in\RR$ and assume that $\varphi \in \dot{W}^{{s,2}}_{\diver}$.   Then $ \overline{\varphi} \in
\dot{W}^{{s}+ \beta,2}_{\diver}$   such that 
\begin{equation}
\label{1suralpharbeta}
\|\overline{\varphi}  \|_{{s}+ \beta,2} \le \frac{1}{\alpha^{ \beta}} \|\varphi \|_{s,2},
  \end{equation}
  and 
  \begin{equation}
\label{1suralphasansbeta}
\|\overline{\varphi}  \|_{{s},2} \le  \|\varphi \|_{s,2}.
  \end{equation}
\end{Lem}
\textbf{Proof.} see in \cite{A01}
  \begin{Lem} 
  Assume  $\vec{w}_{\alpha}$ belongs to $  $ the energy space of  solutions of the Navier-Stokes equations, then   
  \begin{equation}
 \int_{0}^{T}\|{\bw_{\alpha} \otimes \bw_\alpha} \|_{{\frac{8-3r}{2r}},2}^{\frac{r}{2}} dt < \infty \hbox{ for any } \frac{8}{3} \le r < \infty.  
   \end{equation}
   \end{Lem}
   \textbf{Proof.} We have by interpolation that  
    \begin{equation}
   \bw_{\alpha} \in L^{r}(0,T; L^{\frac{6r}{3r-4}}(\tore)^3)
   \end{equation}
    for any $r \ge 2$, thus we deduce by using  H\"{o}lder inequality   that 
    \begin{equation}
    \label{interpolation2}
     {\bw_{\alpha} \otimes \bw_\alpha} \in  L^{\frac{r}{2}}(0,T; L^{\frac{3r}{3r-4}}(\tore)^{3\times 3}).\end{equation}
     From Sobolev embedding we deduce that \begin{equation}
     {\bw_{\alpha} \otimes \bw_\alpha} \in  L^{\frac{r}{2}}(0,T; W^{\frac{8-3r}{2r},2}(\tore)^{3\times 3}),\end{equation}
     for any $r\ge {\frac{8}{3}}.$ 
   \begin{Lem}
   \label{convergncelp}
Assume  $\vec{w}_{\alpha}$ belongs to $  $ the energy space of  solutions of the Navier-Stokes equations, then  for all   $ p \ge 1$ and $q \ge \frac{4}{3} $ such that 
\begin{align}
\label{conditiononpandq}
\frac{1}{p}+\frac{2}{3q} <1,  
\end{align}
we have 

 \begin{equation}
 \label{interpolation1}
\int_{0}^{T}\| \overline{\bw_\alpha \otimes \bw_{\alpha}} -  {\bw_\alpha \otimes \bw_{\alpha}}\|_{p}^{q} dt \le C \alpha^{\frac{3q+p-3pq}{p}}.
 \end{equation}

%\begin{align}
%\label{interpolation1}
%\overline{\bw_\alpha\otimes\bw_\alpha}  &\rightarrow \bw \otimes \bw &&\textrm{ strongly in  }
%L^q(0,T;L^p(\tore)^{3 \times 3})  
%\end{align}
\end{Lem}
\textbf{Proof.} 
We take $r= 2q$, from the Sobolev injection $ \displaystyle {W}^{\frac{3p-6}{2p},2}(\tore)\hookrightarrow L^p(\tore),$ it is sufficent to show that 
  \begin{equation}
\int_{0}^{T}\| \overline{\bw_\alpha \otimes \bw_{\alpha}} -  \bw_{\alpha}\otimes\bw_{\alpha}\|_{{\frac{3p-6}{2p}},2}^{\frac{r}{2}} dt \le C \alpha^{\frac{3r+2p-3pr}{2p}} . 
  \end{equation}
  From the relation between $\overline{\bw_\alpha \otimes \bw_{\alpha}}$  and ${\bw_\alpha \otimes \bw_{\alpha}}$  we have 
    \begin{equation}
  \| \overline{\bw_\alpha \otimes \bw_{\alpha}} - {\bw_\alpha \otimes \bw_{\alpha}}\|_{{\frac{3p-6}{2p}},2}^{\frac{r}{2}} \le \alpha^{ \theta r} \|\overline{\bw_\alpha \otimes \bw_{\alpha}} \|_{{\frac{3p-6}{2p} + 2\theta},2}^{\frac{r}{2}}
   \end{equation}
Lemma \ref{fourierdiscret} implies that 
 \begin{equation}
\int_{0}^{T}\| \overline{\bw_\alpha \otimes \bw_{\alpha}} -  {\bw_\alpha \otimes \bw_{\alpha}}\|_{{\frac{3p-6}{2p}},2}^{\frac{r}{2}} dt \le \alpha^{\frac{3r+2p-3pr}{2p}} \int_{0}^{T}\|{\bw_\alpha \otimes \bw_{\alpha}} \|_{{\frac{8-3r}{2r}},2}^{\frac{r}{2}} dt.
 \end{equation}
Recall that   \begin{equation}
\int_{0}^{T}\|{\bw_{\alpha} \otimes \bw_\alpha} \|_{{\frac{8-3r}{2r}},2}^{\frac{r}{2}} dt < \infty \hbox{ for any } \frac{8}{3} \le r < \infty.  
   \end{equation}
This yields  the desired result for any $ p \ge 1$, $ q \ge \frac{4}{3}$ such that $ \frac{1}{p}+\frac{2}{3q} <1  $.\\
% Whereupon we obtain that 
%\begin{align}
%\overline{\bw_\alpha\otimes\bw_\alpha}  &\rightarrow \bw \otimes \bw &&\textrm{ strongly in  }
%L^{p}(0,T;L^p(\tore)^{3 \times 3}) \textrm{ for all }  \frac{4}{3} \le p  < \frac{5}{3}.
%\end{align}
%   Finally, from (\ref{interpolation2}) and (\ref{interpolation1})  we  deduce   that 
%\begin{align}
%\overline{\bw_\alpha\otimes\bw_\alpha}  &\rightharpoonup \bw \otimes \bw &&\textrm{ weakly in  }  L^{\frac{r}{2}}(0,T; L^{\frac{3r}{3r-4}}(\tore)^{3\times 3}) \textrm{ for all }   r  \ge 2.
%\end{align}

%\BEQ  \label{GQAP954} 
% || \overline{\bu} ||_{s+2\theta} \le {1 \over \alpha^{2 \theta}} || \bu ||_{\vec{H}^{s}}.  \EEQ 
% We shall sometimes denote $\overline{\bu_\alpha}$ instead of $\overline{\bu}$, when we need to recall the dependance on the $\alpha$ parameter. 

\textbf{Proof of Theorem \ref{1deuxieme}.}
The proof of Theorem \ref{1deuxieme} follows the lines of the proof of the Theorem 4 in \cite{LL06b}.
The only difference is the strong convergence of the pressure term $q_{\alpha}$ to the pressure term $p$ of the Navier-Stokes equations. 
 We will use Layton-Lewandowski  \cite{LL06b} as a
reference and only point out the differences between their  proof of convergence to a weak
solution of the Navier-Stokes equations and the proof of convergence  in our study.  
First, we need to find estimates
that are independent from $\alpha$. 
%De la relation entre $\bu_\alpha$ et $\overline{\bu_\alpha}$ on peut montrer que pour tout $\bu_\alpha \in L^p(0,T;L^p(\tore)^3) $ lorsque $\alpha$ tend vers 0 on a :  
%\begin{align}
%\overline{\bu_\alpha}  &\rightarrow \bu &&\textrm{strongly in  }
%L^p(0,T;L^p(\tore)^3) \textrm{ pour tout } p>1.
%\end{align}
%(Voir Lemme 3.1 page 8)
Using the fact that  $\vec{w}_{\alpha}$  belong to the energy space:  $L^{\infty}
([0,T];L^2({\tore})^3)\cap L^2([0,T];\dot{W}^{1,2}_{\diver})$ and from the  Aubin-Lions compactness Lemma  (the same arguments as in section 2.1) we can find a  subsequence $(\vec{w}_{\alpha_j}, q_{\alpha_j})$  and $(\vec{v}_{}, p_{})$ such that when   $\alpha_j$ tends to zero we have: 
%La solution $(\vec{u}_{\alpha_j}, p_{\alpha_j})$  appartient a l'espace d'énergie de Leray donc d'aprés le Lemme de compacité de Aubin-Lions on peut extraire une sous suite tel que lorque $\alpha_j$ tend vers 0 on a: 
\begin{align}
\bw_{\alpha_j} &\rightharpoonup^* \bv &&\textrm{weakly$^*$ in } L^{\infty}
([0,T];L^2(\tore)^3), \label{hanic122}\\
%\overline{\bv}^\epsilon &\rightharpoonup^* \overline{\bv} &&\textrm{weakly$^*$ in } L^{\infty}
%(0,T;W^{\frac{1}{2},2}_{}), \label{c122prime}\\
\bw_{\alpha_j} &\rightharpoonup \bv &&\textrm{weakly in }
L^2([0,T];W^{1,2}(\tore)^3), \label{hanic22}\\
\bw_{\alpha_j} &\rightharpoonup \bv &&\textrm{weakly in }  L^{r}(0,T; L^{\frac{6r}{3r-4}}(\tore)^3) \textrm { for all } r \ge 2, \label{hanic22rrrr}\\
%\overline{\bv}^\epsilon &\rightharpoonup \overline{\bv} &&\textrm{weakly in }
%L^2(0,T;W^{\frac{3}{2},2}_{})\cap L^{5}(0,T;L^{5}(\tore)^3), \label{c22prime}\\
%\bv^\epsilon_{,t}&\rightharpoonup \bv_{,t} &&\textrm{weakly in } L^{2}
%(0,T;W^{-1,2}_{}),
%\label{c322}\\
{\bw_{\alpha_j}\otimes\bw_{\alpha_j}}  &\rightharpoonup \bv\otimes\bv &&\textrm{weakly in  }
 L^{\frac{r}{2}}(0,T; L^{\frac{3r}{3r-4}}(\tore)^{3 \times 3}), \textrm { for all } r \ge 2,\\
%q_{\alpha_j}&\rightharpoonup p &&\textrm{weakly in }  L^{\frac{r}{2}}(0,T; L^{\frac{3r}{3r-4}}(\tore)), \textrm { for all } r \ge 2 \label{hanic32}\\
%\bv^\epsilon &\rightharpoonup \bv &&\textrm{weakly in } L^{\frac{8}{3}}
%(0,T;L^{\frac{8}{3}}(\partial \Omega)^3). \label{c5.12}\\
\bw_{\alpha_j} &\rightarrow \bv &&\textrm{strongly in }
L^p([0,T];L^p(\tore)^3) \textrm{ for all } 2 \le  p<\frac{10}{3},\label{hanic82}\\
%\overline{\bv}^\epsilon &\rightarrow \overline{\bv} &&\textrm{strongly in  }
%L^2(0,T;L^6(\tore)^3),\label{c82pppp}
%\overline{\bv}^\epsilon &\rightarrow \overline{\bv} &&\textrm{strongly in  }
%L^q(0,T;L^q(\tore)^3) \textrm{  for all } q< 5 .\label{c82pppp}
%\bv^\epsilon &\rightarrow \bv &&\textrm{strongly in  } L^q(0,T;L^q(\partial
%\Omega)^3) \textrm{ for all } q<\frac{8}{3},\label{c5.22}
{\bw_{\alpha_j}\otimes\bw_{\alpha_j}} &\rightarrow \bv\otimes\bv &&\textrm{strongly in  }
L^p([0,T];L^p(\tore)^{3\times 3}) \textrm{ for all } \frac{4}{3} \le p < \frac{5}{3},\label{pkha}
%q_{\alpha_j}&\rightarrow p &&\textrm{strongly in  }
%L^p([0,T];L^p(\tore)) \textrm{ for all }  \frac{4}{3}  \le p < \frac{5}{3}.
\end{align}
Having (\ref{interpolation1}) and (\ref{pkha}) at hand we deduce that 
%From lemma \ref{} and we deduce that for any $p$ , $q$ verifies \ref{conditiononpandq} 
%\begin{align}
%\label{interpolation1}
%\overline{\bw_\alpha\otimes\bw_\alpha}  &\rightarrow \bv \otimes \bv &&\textrm{ strongly in  }
%L^q(0,T;L^p(\tore)^{3 \times 3})  
%\end{align}
%Whereupon we obtain that 
\begin{align}\label{interpolationpkha}
\overline{\bw_{\alpha_j}\otimes\bw_{\alpha_j}}  &\rightarrow \bv \otimes \bv &&\textrm{ strongly in  }
L^{p}(0,T;L^p(\tore)^{3 \times 3}) \textrm{ for all }  \frac{4}{3} \le p  < \frac{5}{3}.
\end{align}
   Then, from (\ref{interpolation2}) and (\ref{interpolationpkha})  we  deduce   that 
\begin{align}
\label{hanimodofie}
\overline{\bw_{\alpha_j}\otimes\bw_{\alpha_j}}  &\rightharpoonup \bv \otimes \bv &&\textrm{ weakly in  }  L^{\frac{r}{2}}(0,T; L^{\frac{3r}{3r-4}}(\tore)^{3\times 3}) \textrm{ for all }   r  \ge 2.
\end{align}
Further, we have 
\begin{align}
q_{}(t)=\vec{R}(\sum_{k,l}\overline{\bw^k_{{\alpha_j}} \bw^l_{\alpha_j}}(t) )
\end{align}
where the linear map $\vec{R}$ is defined by 
\begin{align}
\vec{R} \quad  : \quad &L^s(\tore)^9 \longmapsto L^s(\tore)\\
    & (u^{kl})_{k,l=1,2,3}\longmapsto  (-\Delta)^{-1}\partial_k \partial_l (u^{kl}) 
\end{align}
 By the theory of Riesz transforms, $\vec{R}$ is a continuous map for any $s \in ]1, \infty[$. 
  Consequently, from (\ref{hanimodofie}) we have  
\begin{align}
\label{hihihihi}
\int_{0}^{T}\| q_{\alpha_j}\|_{\frac{3r}{3r-4}}^{\frac{r}{2}}dt < \infty, \textrm{ for all }   r  \ge 2.
%q_{\alpha_j}&\rightharpoonup p &&\textrm{weakly in }  L^{\frac{r}{2}}(0,T; L^{\frac{3r}{3r-4}}(\tore)), \textrm { for all } r \ge 2 \label{hanic32}
\end{align}
From (\ref{interpolationpkha}) we deduce that  for almost all $t > 0$, 
\begin{align}
\overline{\bw_{\alpha_j}\otimes\bw_{\alpha_j}}(t)  &\rightarrow \bv \otimes \bv (t) &&\textrm{ strongly in  }
L^p(\tore)^{3 \times 3} \textrm{ for all }  \frac{4}{3} \le p  < \frac{5}{3}.
\end{align}
Using the dominate convergence theorem  and the continuity
of the operator $\vec{R}$, we conclude that 
\begin{align}
q_{\alpha_j}&\rightarrow p &&\textrm{strongly in  }
L^p([0,T];L^p(\tore)) \textrm{ for all }  \frac{4}{3}  \le p < \frac{5}{3}. \label{hhahahaha}
\end{align}
Finally we deduce from (\ref{hhahahaha}) and (\ref{hihihihi}) that 
\begin{align}
q_{\alpha_j}&\rightharpoonup p &&\textrm{weakly in }  L^{\frac{r}{2}}(0,T; L^{\frac{3r}{3r-4}}(\tore)), \textrm { for all } r \ge 2. \label{hanic32}
\end{align}

These convergence results allow us  to prove in the same way as in \cite{LL06b}  that $ (\bv,p)$ is a weak solution to the Navier-Stokes equations, so we will not repeat it.

%\subsection{Consistency error in the critical case } 
%Next, we turn to the consistency errors in the critical case $\theta=\frac{1}{6}$. We will show that the critical regularization leads to a sharper bounds on the consistency error. 
%Let $\tau_{\bv} := \overline{\bv}^{\frac{1}{6}} \otimes \overline{\bv}^{\frac{1}{6}} - \bv \otimes \bv  $  denotes the model  consistency error,
%where $\bv$ is the velocity part of the  solution of the NSE equations obtained as a limit of a subsequence of the sequence $\bw$. We prove that 
%that $\|\overline{\bv}^{\frac{1}{6}}- \bw \|_{L^{\infty}(0,T; L^{2})}$  are bounded by $\|\tau_{\bv}\|_{L^{2}(0,T; L^{2})}$.

\section{Application to the LES  for magnetohydrodynamic equations (LES for MHD)}

In this section, we consider the critical LES regularization for  magnetohydrodynamic (LES for MHD) equations, given by
\begin{align}
 \partial_t\bw-\nu_1 \Delta \bw +\diver(\overline{\bw\otimes\bw}^{\frac{1}{6}})-\diver(\overline{\vec{W}\otimes\vec{W}}^{\frac{1}{6}})
+\nabla q&=0,
\label{2EMEIMHDV1}\\
\partial_t \vec{W} -\nu_2 \Delta \vec{W} +\diver(\overline{\bw\otimes\vec{W}}^{\frac{1}{6}})-\diver(\overline{\vec{W}\otimes\bw}^{\frac{1}{6}})&=0,
\label{2EMEIMHDV2}\\
\displaystyle\int_{\tore} \vec{W}\,d\bx =  \int_{\tore} \bw\,d\bx =0, \quad \diver \vec{w}= \diver \vec{W} &= 0,\label{2EMEIMHDV3}\\
 \vec{W}(0)=\vec{W}_{0},\; \bw(0)&=\bw_{0},\label{2EMEIMHDV4}
\end{align}
where the boundary conditions are taken to be periodic,   and we take  as before the same spacing average operator.
\begin{align}
\alpha^{\frac{1}{3}}( -\Delta)^{\frac{1}{6}} \overline{\varphi}^{\frac{1}{6}} +   \overline{\varphi}^{\frac{1}{6}}
&=\varphi, \quad \varphi(t, \bx + L \vec{e}_{\vec{j}})=\varphi(t,\bx) .\label{TKE}
\end{align}
 Here, the unknowns are the averaging fluid velocity field $\bw(t,\vec{x})$, the averaging fluid pressure $q(t,\vec{x})$, and the averaging  magnetic field $\vec{W}(t,\vec{x})$.  Note that when $\alpha=0$, we formally retrieve the  MHD equations.
Existence and uniqueness results for MHD equations are established by G. Duvaut and J.L. Lions in \cite{DuLi72}. These results are completed by M. Sermange and R. Temam in \cite{RT83mhd}. They showed that the classical properities of the Navier-Stokes equations can be extended to the MHD system. 

The aim in this section is to extend the results of existence uniqueness and convergence established above for the  LES for NSE to the LES for MHD.  We know, thanks  to  the work \cite{latr2010}, that for $\theta=1$  these results hold ture.  Further, 
when $\theta=\frac{1}{6}$, we  proved in the above section  the existence  of a unique  ``regular'' weak solution to the LES for NSE.
Therefore, it is intersecting  to find the critical value of   regularization  needed to establish global in time existence of a unique ``regular'' weak solution to LES for MHD.

We divide this section into two subsections.  One is devoted to prove the  existence  of a unique  ``regular'' weak solution to the LES for MHD with $\theta=\frac{1}{6}$.  The second one is devoted to prove that this solution converges to a weak solution to the MHD equations when $\alpha$ tends to zero. 
%%The ladder inequality for the MHD equation are showed in \cite{}.\\
% The use of Leray-$\alpha$ regularization to the MHD equations has received many studies see \cite{LT2007}. The idea to consider the LES for MHD with critical regularization  is a new feature for the present work.\\
%%Our present study is largely inspired by all this work. 
%In this approach
 %we will establesh  the global existence and uniqueness of solutions for the MHD-Deconvolution equations (\ref{IMHDV}) for $\theta=1/4$.\\ 
\subsection{Existence and uniqueness results for the LES for MHD }
%\subsection{ Global existence and uniqness for the MHD-Deconvolution equations}
First, we establish  the global existence and uniqueness of solutions for the LES for MHD equations with  $\theta=\frac{1}{6}$.\\
We have the following theorem:
\begin{Theorem}
\label{existensssIMHDV}
 Assume that $\theta = \frac{1}{6}$. Assume $\bw_0$ and  $\bW_0$ are both in  $ W^{\frac{1}{6},2}_{\diver}.$  Then there exist $(\bw,\vec{W} ,q)$ a unique ``regular"  weak  solution to \eqref{2EMEIMHDV1}--\eqref{2EMEIMHDV4}  such that
\begin{align}
\bw, \ \vec{W} &\in \mathcal{C}_{}(0,T;\dot{W}^{\frac{1}{6},2}_{\diver}) \cap
L^2(0,T;\dot{W}^{1+\frac{1}{6},2}_{\diver}),\label{mhdbv12}\\
\bw_{,t}, \  \vec{W}_{,t}&\in   L^{2}(0,T;W^{-\frac{5}{6},2}_{}),
\label{mhdbvt}\\
q&\in L^{2}(0,T;W^{\frac{1}{6},2}(\tore))
\label{mhdpsp}.
\end{align}
fulfill
\begin{equation}
\begin{split}
\int_0^T \langle \bw_{,t}, \bfi \rangle  -  (\overline{\bw \otimes \bw}, \nabla
\bfi)+ (\overline{\bW \otimes \bW}, \nabla
\bfi) +  \nu_1(
\nabla \bw, \nabla \bfi  )\; dt =0\\
\int_0^T \langle \bW_{,t}, \bfi \rangle  -  (\overline{\bw \otimes \bW}, \nabla
\bfi)+ (\overline{\bW \otimes \bw}, \nabla
\bfi) +  \nu_2(
\nabla \bw, \nabla \bfi  )\; dt =0\\
\qquad \textrm{ for all } \bfi\in L^{2}(0,T; W^{\frac{5}{6},2}_{\diver}).
\end{split}\label{mhdweak1}
\end{equation}
Moreover,
\begin{equation}
\bw(0)= \bw_0 \quad  and \quad \bW(0)=\bW_0. 
\label{mhdintiale}
\end{equation}
 %and the initial conditions are attained in the following sense
%\begin{equation}
%\lim_{t\to 0+}\left(\|\bw(t)-\bw_0\|_2^2 \right)
%=0. \label{inca}
%\end{equation}
\end{Theorem}
 \textbf{Proof. of Theorem \ref{existensssIMHDV}.} 
 We only sketch  the proof since  is similar to the   Navier-Stokes equations case. The proof  is obtained by taking the inner product of (\ref{IMHDV1}) with $\alpha^{\frac{1}{3}}(-\Delta)^{\frac{1}{6}}\vec{w}+ \vec{w},$ (\ref{IMHDV2}) with $\alpha^{\frac{1}{3}}(-\Delta)^{\frac{1}{6}}\vec{W}+ \vec{W}$ and then adding them, the existence of a
solution to the critical LES for MHD  can be derived thanks to the Galerkin
method. Notice that  $({\vec{w}^{}}, {\bW^{ }})$    satisfy the following estimates
\begin{equation}
\begin{array}{lllll}
 \label{mhdapriori12leray}
\displaystyle \left(\|{\bw} (t)\|_{2}^2 + \|{\bw}(t)\|_{\frac{1}{6},2}^2\right)+\left(\|{\bW} (t)\|_{2}^2
+ \|{\bW}(t)\|_{\frac{1}{6},2}^2\right)\\
\quad +\displaystyle 2\nu_1\int_{0}^{t}\left(\|{\bw} \|_{1,2}^2  + \|{\bw} \|_{1+\frac{1}{6},2}^2 \right)ds +2\nu_2\int_{0}^{t}\left(\|{\bW} \|_{1,2}^2  + \|{\bW} \|_{1+\frac{1}{6},2}^2 \right)ds \\
\quad \quad \quad = \left(\|{\bw}_0 \|_{2}^2 + \|{\bw}_0\|_{\frac{1}{6},2}^2\right)+\left(\|{\bW}_0 \|_{2}^2 + \|{\bW}_0\|_{\frac{1}{6},2}^2\right).
\end{array}
\end{equation}
The averaging pressure $q^{ }$   is reconstructed from $\vec{w}$ and $\bW $  (as we work with periodic boundary conditions) 
%thanks the De Rham Theorem 
and its regularity
results from the fact that  $ \overline{\bw \otimes \bw}$ and    $\overline{\bW \otimes \bW} \in L^2([0,T];{W}^{\frac{1}{6},2}(\tore)^{3\times 3})$.\\
It remains to prove the uniqueness. Let $(\vec{w}_{1},\vec{W}_{1}, q_{1})$ and $(\vec{w}_{2},\vec{W}_{2},
q_{2})$,  be two solutions, $\delta\vec{w} = \vec{w}_{2}-\vec{w}_{1}$,$\delta\bW = \bW_{2}-\bW_{1}$, $\delta q =
q_{2}-q_{1}$. Then one has
\begin{equation}
\label{matin1}
\begin{split}
\partial_{t} \delta\vec{w}- \nu_1
\Delta\delta\vec{w} + \diver(\overline{\bw_2 \otimes \bw_2})- \diver(\overline{\bw_1 \otimes \bw_1})\\
\quad \quad  - \diver(\overline{\bW_2 \otimes \bW_2})+\diver(\overline{\bW_1 \otimes \bW_1})  + \nabla\delta q =0,\\
\partial_{t} \delta\bW  - \nu_2
\Delta\delta\bW + \diver(\overline{\bw_2 \otimes \bW_2})- \diver(\overline{\bw_1 \otimes \bW_1})\\
\quad \quad  - \diver(\overline{\bW_2 \otimes \bw_2})+\diver(\overline{\bW_1 \otimes \bw_1})    = 0,
\end{split}
\end{equation}
and $\delta\vec{w} = 0$, $\delta\vec{W} =0$  at initial time. One can take $\alpha^{\frac{1}{3}}(-\Delta)^{\frac{1}{6}}\delta\vec{w}+ \delta\vec{w}$
 as test in  the first equation of (\ref{matin1}) and $\alpha^{\frac{1}{3}}(-\Delta)^{\frac{1}{6}}\delta\vec{W}+ \delta\vec{W}$
 as test in the second equations of (\ref{matin1}). Since $
\vec{w}_1$ is divergence-free we have 
\begin{equation}
\begin{split}\displaystyle
\int_{0}^{T} \int_{\tore}   \vec{w}_1 \otimes \delta \vec{w} : \nabla \delta\vec{w}= -\int_{0}^{T} \int_{\tore}   (\vec{w}_1 \cdot \nabla)\delta \vec{w} \cdot \delta\vec{w} =0  ,\\
 \end{split}
\end{equation}
Thus we obtain by using the fact that the averaging operator  commutes with differentiation under periodic boundary conditions 
 \begin{equation}
 \label{soira}
\begin{split}\displaystyle
\int_{0}^{T} \int_{\tore}\left(\diver(\overline{\bw_2 \otimes \bw_2})-\diver(\overline{\bw_1 \otimes \bw_1})\right) \cdot\left( \alpha^{\frac{1}{3}}(-\Delta)^{\frac{1}{6}}\delta\vec{w}+ \delta\vec{w} \right)\\
= \int_{0}^{T} \int_{\tore}\left(\diver({\bw_2 \otimes \bw_2})-\diver({\bw_1 \otimes \bw_1})\right) \cdot\delta\vec{w}\\ 
= -\int_{0}^{T} \int_{\tore}  \delta \vec{w} \otimes \vec{w}_2 : \nabla \delta\vec{w}.
 \end{split}
\end{equation}
Similarly, 
because 
$
(\vec{w}_1)$ is divergence-free we have 
\begin{equation}
\begin{split}\displaystyle
 \int_{0}^{T} \int_{\tore}   \vec{w}_1 \otimes \delta \vec{W} : \nabla \delta\vec{W}= -\int_{0}^{T} \int_{\tore}   (\vec{w}_1 \cdot \nabla) \delta \vec{W} \cdot \delta\vec{W}  =0  ,\\
 \end{split}
\end{equation}
and thus we have the following identity 
 \begin{equation}
\begin{split}\displaystyle
\int_{0}^{T} \int_{\tore}\left(\diver(\overline{\bw_2 \otimes \bW_2})-\diver(\overline{\bw_1 \otimes \bW_1})\right) \cdot\left( \alpha^{\frac{1}{3}}(-\Delta)^{\frac{1}{6}}\delta\vec{W}+ \delta\vec{W} \right)\\
= \int_{0}^{T} \int_{\tore}\left(\diver({\bw_2 \otimes \bW_2})-\diver({\bw_1 \otimes \bW_1})\right) \cdot\delta\vec{W}\\ 
= -\int_{0}^{T} \int_{\tore}  \delta \vec{w} \otimes \vec{W}_2 : \nabla \delta\vec{W}.\\
 \end{split}
\end{equation}
Concerning the remaining terms 
we get by integrations by parts and by using the  using the fact that the averaging operator commutes with differentiation under periodic boundary conditions
\begin{equation}
\begin{split}\displaystyle
\int_{0}^{T} \int_{\tore}\left(-\diver(\overline{\bW_2 \otimes \bW_2})+\diver(\overline{\bW_1 \otimes \bW_1})\right) \cdot\left( \alpha^{\frac{1}{3}}(-\Delta)^{\frac{1}{6}}\delta\vec{w}+ \delta\vec{w} \right)\\
= \int_{0}^{T} \int_{\tore}\left(-\diver({\bW_2 \otimes \bW_2})+\diver({\bW_1 \otimes \bW_1})\right) \cdot\delta\vec{w}\\ 
= \int_{0}^{T} \int_{\tore}   \vec{W}_1 \otimes \delta \vec{W} : \nabla \delta\vec{w} + \int_{0}^{T} \int_{\tore}  \delta \vec{W} \otimes \vec{W}_2 : \nabla \delta\vec{w}.\\
 \end{split}
\end{equation}
and similarly 
\begin{equation}
\label{soirb}
\begin{split}\displaystyle
\int_{0}^{T} \int_{\tore}\left(-\diver(\overline{\bW_2 \otimes \bw_2})+\diver(\overline{\bW_1 \otimes \bw_1})\right) \cdot\left( \alpha^{\frac{1}{3}}(-\Delta)^{\frac{1}{6}}\delta\vec{W}+ \delta\vec{W} \right)\\
= \int_{0}^{T} \int_{\tore}\left(\diver({\bW_2 \otimes \bw_2})-\diver({\bW_1 \otimes \bw_1})\right) \cdot\delta\vec{W}\\ 
= -\int_{0}^{T} \int_{\tore}   (\vec{W}_1 \cdot \nabla) \delta \vec{w} \cdot  \delta\vec{W} + \int_{0}^{T} \int_{\tore}  \delta \vec{W} \otimes \vec{w}_2 : \nabla \delta\vec{W}.\\
 \end{split}
\end{equation}

%
%Therefore, 
%
% \begin{equation}
%\label{soira}\begin{array}{ll}
% \displaystyle {\frac{d }{2dt}} \int_{\tore} | \delta\vec{u}| ^{2} +\displaystyle \nu_1\int_{\tore}| \nabla
%\delta\vec{u} |^{2}  -\displaystyle \int_{\tore} (H_{N} (\B_{1}) \nabla) \delta\B. \delta
%\vec{u}\\
% \hskip 3cm  =\displaystyle -\int_{\tore}(H_{N} (\delta\vec{u}) \nabla) \vec{u}_{2}. \delta
%\vec{u} +\displaystyle \displaystyle\int_{\tore}(H_{N} (\delta\B) \nabla) \B_{2}. \delta
%\vec{u},
%\end{array}
%\end{equation}
%and 
%  \begin{equation}
%\label{soirb}\begin{array}{ll}
% \displaystyle {\frac{d }{2dt}} \int_{\tore}  | \delta\B| ^{2} + \displaystyle\nu_2\int_{\tore} | \nabla
%\delta\B |^{2}   - \displaystyle \int_{\tore} (H_{N} (\B_{1}) \nabla) \delta\vec{u} . \delta
%\B\\ 
% \hskip 3cm  = \displaystyle \int_{\tore}(H_{N} (\delta\B) \nabla) \vec{u}_{2}     . \delta
%\B   \displaystyle   -\displaystyle\int_{\tore}   (H_{N} (\delta\vec{u}) \nabla) \B_{2}  . \delta
%\B.
%\end{array}
%\end{equation}
%
%One has by a integration by parts,
%\BEQ \label{MATIN2}
%- \displaystyle \int_{\tore} (H_{N} (\B_{1}) \nabla) \delta\B. \delta
%\vec{u}     = \displaystyle \int_{\tore} (H_{N} (\B_{1}) \nabla) \delta\vec{u} . \delta
%\B. 
%\EEQ
Therefore by adding (\ref{soira})-(\ref{soirb}) and using  the fact that the averaging operator commutes with differentiation under periodic boundary conditions we obtain
\begin{equation}
 \label{soirab}
\begin{array}{llll}
 \displaystyle {\frac{d }{2dt}} \int_{\tore} \left(| \delta\vec{w}| ^{2} + \alpha^{\frac{1}{6}}| \nabla^{\frac{1}{6}}\delta\vec{w}| ^{2}  \right)+{\frac{d }{2dt}}\left( \int_{\tore}  | \delta\bW| ^{2}  +\alpha^{\frac{1}{6}}| \nabla^{\frac{1}{6}}\delta\vec{W}| ^{2}\right)\\ +\displaystyle \nu_1\left(\int_{\tore}| \nabla
\delta\vec{u} |^{2} + | \nabla^{1+\frac{1}{6}}
\delta\vec{u} |^{2} \right)+\displaystyle\nu_2\left(\int_{\tore} | \nabla
\delta\B |^{2}+| \nabla^{1+\frac{1}{6}}
\delta\B |^{2} \right)\\
\hskip 2cm=\displaystyle \int_{\tore}\delta \vec{w} \otimes \vec{w}_2 : \nabla \delta\vec{w}   + \displaystyle \displaystyle \int_{\tore}  \delta \vec{w} \otimes \vec{W}_2 : \nabla \delta\vec{W}  \\
\hskip 4cm
- \displaystyle \int_{\tore}   \delta \vec{W} \otimes  \vec{W}_2 : \nabla \delta\vec{w}  - \displaystyle  \int_{\tore}  \delta \vec{W} \otimes  \vec{w}_2 : \nabla \delta\vec{W}.
 \end{array}
\end{equation}

%
%
%One has by a  integration by parts,
%\begin{subequations}
%\begin{align}
%-\int_{\tore}(H_{N} (\delta\vec{u}) \nabla) \vec{u}_{2}. \delta
%\vec{u} =\int_{\tore} H_{N} (\delta\vec{u}) \otimes \vec{u}_{2} : \g \delta
%\vec{u},\\
%\int_{\tore}(H_{N} (\delta\B) \nabla) \B_{2}. \delta
%\vec{u} = -\int_{\tore} H_{N} (\delta\B) \otimes \B_{2} : \g \delta
%\vec{u},\\
%\int_{\tore}(H_{N} (\delta\B) \nabla) \vec{u}_{2}. \delta
%\B = -\int_{\tore} H_{N} (\delta\B) \otimes \vec{u}_{2} : \g \delta
%\B,\\
%-\int_{\tore}(H_{N} (\delta\vec{u}) \nabla) \B_{2}. \delta
%\B = \int_{\tore} H_{N} (\delta\vec{u}) \otimes \B_{2} : \g \delta
%\B.
%\end{align}
%\end{subequations}
By the norm duality
\begin{align}
|\int_{\tore}\delta \vec{w} \otimes \vec{w}_2 : \nabla \delta\vec{w} |  \le    \| \delta \vec{w} \otimes \vec{w}_2  \|_{-\frac{1}{6},2}\|\nabla \delta\vec{w} \|_{\frac{1}{6},2},\\
 |\int_{\tore}    \delta \vec{w} \otimes \vec{W}_2 : \nabla \delta\vec{W} |  \le   \| \delta \vec{w} \otimes \vec{W}_2  \|_{-\frac{1}{6},2}\|\nabla \delta\vec{W} \|_{\frac{1}{6},2},   \\
|\int_{\tore}    \delta \vec{W} \otimes  \vec{W}_2 : \nabla \delta\vec{w}      |  \le  \|  \delta \vec{W} \otimes  \vec{W}_2  \|_{-\frac{1}{6},2}\|\nabla \delta\vec{w} \|_{\frac{1}{6},2},  \\
|\int_{\tore}     \delta \vec{W} \otimes  \vec{w}_2 : \nabla \delta\vec{W}     |  \le     \|   \delta \vec{W} \otimes  \vec{w}_2  \|_{-\frac{1}{6},2}\|\nabla \delta\vec{W} \|_{\frac{1}{6},2}.
 \end{align}

By Young's inequality,

\begin{align}
|\int_{\tore}\delta \vec{w} \otimes \vec{w}_2 : \nabla \delta\vec{w} |  \le    \frac{1}{\nu_1}\| \delta \vec{w} \otimes \vec{w}_2  \|_{-\frac{1}{6},2}^2 + \frac{\nu_1}{4}\|\nabla \delta\vec{w} \|_{\frac{1}{6},2}^2 ,\\
 |\int_{\tore}    \delta \vec{w} \otimes \vec{W}_2 : \nabla \delta\vec{W} |  \le  \frac{1}{\nu_2} \| \delta \vec{w} \otimes \vec{W}_2  \|_{-\frac{1}{6},2}^2+ \frac{\nu_2}{4}\|\nabla \delta\vec{W} \|_{\frac{1}{6},2}^2,   \\
|\int_{\tore}    \delta \vec{W} \otimes  \vec{W}_2 : \nabla \delta\vec{w}      |  \le  \frac{1}{\nu_1}\|  \delta \vec{W} \otimes  \vec{W}_2  \|_{-\frac{1}{6},2}^2 +\frac{\nu_1}{4}\|\nabla \delta\vec{w} \|_{\frac{1}{6},2}^2,  \\
|\int_{\tore}     \delta \vec{W} \otimes  \vec{w}_2 : \nabla \delta\vec{W}     |  \le    \frac{1}{\nu_2} \|   \delta \vec{W} \otimes  \vec{w}_2  \|_{-\frac{1}{6},2}^2+  \frac{\nu_2}{4}\|\nabla \delta\vec{W} \|_{\frac{1}{6},2}^2.
 \end{align}

%\begin{align}
%|\int_{\tore}\delta \vec{w} \otimes \vec{w}_2 : \nabla \delta\vec{w} |  \le \frac{\nu_1}{4} \int_{\tore} |\nabla \delta\vec{w} |^2 + \frac{1}{\nu_1} \int _{\tore}| \delta \vec{w} \otimes \vec{w}_2  |^2,\\
% |\int_{\tore}(H_{N} (\delta\B) \nabla) \B_{2}. \delta
%\vec{u}|  \le \frac{\nu_1}{4} \int_{\tore} |\g \delta \vec{u} |^2 + \frac{1}{\nu_1} \int _{\tore}| H_N (\delta \B) |^2 | \B_2|^2,\\
%|\int_{\tore}(H_{N} (\delta\B) \nabla) \vec{u}_{2}. \delta
%\B|  \le \frac{\nu_2}{4} \int_{\tore} |\g \delta \B |^2 + \frac{1}{\nu_2} \int _{\tore}| H_N (\delta \B) |^2 | \vec{u}_2|^2,\\
%|\int_{\tore}(H_{N} (\delta\vec{u}) \nabla) \B_{2}. \delta
%\B|  \le \frac{\nu_2}{4} \int_{\tore} |\g \delta \B |^2 + \frac{1}{\nu_2} \int _{\tore}| H_N (\delta \vec{u}) |^2 | \B_2|^2.
% \end{align}

By H\"{o}lder inequality combined with Sobolev injection

\begin{align}
\begin{array}{llll}
 \displaystyle   \frac{1}{\nu_1}\| \delta \vec{w} \otimes \vec{w}_2  \|_{-\frac{1}{6},2}^2    &\le \displaystyle \frac{1}{\nu_1} \| \delta \vec{w} \|_{\frac{1}{6},2}^2 \|  \vec{w}_2 \|_{1+\frac{1}{6},2}^2&\\
 \end{array} \\
\begin{array}{llll}
 \frac{1}{\nu_2} \| \delta \vec{w} \otimes \vec{W}_2  \|_{-\frac{1}{6},2}^2
 &\le \displaystyle  \frac{1}{\nu_2} \| \delta \vec{w} \|_{\frac{1}{6},2}^2 \|  \vec{W}_2 \|_{1+\frac{1}{6},2}^2&\\
 \end{array}\\
\begin{array}{llll}
 \displaystyle \frac{1}{\nu_1}\|  \delta \vec{W} \otimes  \vec{W}_2  \|_{-\frac{1}{6},2}^2 
 &\le \displaystyle \frac{1}{\nu_1} \| \delta \vec{W} \|_{\frac{1}{6},2}^2 \|  \vec{W}_2 \|_{1+\frac{1}{6},2}^2&\\
\end{array}\\
\begin{array}{llll}
 \displaystyle  \frac{1}{\nu_2} \|   \delta \vec{W} \otimes  \vec{w}_2  \|_{-\frac{1}{6},2}^2
 &\le \displaystyle \frac{1}{\nu_2} \| \delta \vec{W} \|_{\frac{1}{6},2}^2 \|  \vec{w}_2 \|_{1+\frac{1}{6},2}^2  &
\end{array}
 \end{align}

Hence,

\begin{equation}
 \label{soirabc}
\begin{array}{llll}
 \displaystyle {\frac{d }{2dt}} \int_{\tore} \left( \| \delta \vec{w} \|_{2}^2  +\displaystyle   \| \delta \bW \|_{2}^2   \right)+ \alpha^{\frac{1}{6}}{\frac{d }{2dt}}\left(  \| \delta \vec{w} \|_{\frac{1}{6},2}^2  +\displaystyle   \| \delta \bW \|_{\frac{1}{6},2}^2    \right)\\ +\displaystyle \min{(\nu_1,\nu_2)} \left(   \| \delta \vec{w} \|_{1,2}^2  +\displaystyle   \| \delta \bW \|_{1,2}^2   \right)+\displaystyle \alpha^{\frac{1}{6}} \min{(\nu_1,\nu_2)}\left(   \| \delta \vec{w} \|_{1+\frac{1}{6},2}^2  +\displaystyle   \| \delta \bW \|_{1+\frac{1}{6},2}^2        \right)\\
\hskip 1cm \displaystyle \le \frac{1}{\min{(\nu_1,\nu_2)}}\left(
 \| \delta \vec{w} \|_{\frac{1}{6},2}^2  +\displaystyle   \| \delta \bW \|_{\frac{1}{6},2}^2 \right)\left(\| \vec{w}_2\|_{1+\frac{1}{6},2}^2 +   \| \bW_2\|_{1+\frac{1}{6},2}^2    \right)
 \end{array}
\end{equation}
We conclude that $\delta\vec{u} =\delta\B= 0$ thanks to Gr\"{o}nwall's Lemma.\\

\subsection{Limit consistensy for the critical LES for 	MHD }
Next, we will deduce that the LES for MHD with critical  regularization   gives rise to a  weak solution to the MHD equations.
 \begin{Theorem}
 \label{conmhd}
 Let $(\bw_{\alpha}, \bW_{\alpha},q_{\alpha})$ be the solution of \eqref{BM}--\eqref{vnzTKE}  for a fixed $\alpha$.
 There is a subsequence $\alpha_j$  such that
$(\bw_{\alpha_j}, \bW_{\alpha_j},q_{\alpha_j}) \rightarrow (\bv, \B, p)$ as $j \rightarrow \infty $
where
$ (\bv, \bW, p)  \in [L^{\infty}
([0,T];L^2({\tore})^3)\cap L^2([0,T];\dot{W}^{1,2}_{\diver})]^2 \times L^{\frac{5}{3}}([0,T];L^{\frac{5}{3}}(\tore))$
is a weak solution of the Navier-Stokes equations  with periodic boundary conditions and zero mean value constraint.\\
 The sequence  $ \bw_{\alpha_j}$  converges strongly
to $\bv$ in the space $L^p([0,T];L^p(\tore)^3)$
for all  $2 \le  p<\frac{10}{3}, $  and weakly in
$ L^{r}(0,T; L^{\frac{6r}{3r-4}}(\tore)^3) \textrm { for all } r \ge 2. $\\
The sequence  $ \bW_{\alpha_j}$  converges strongly
to $\B$ in the space $L^p([0,T];L^p(\tore)^3)$ for all  $2 \le  p<\frac{10}{3}, $  and weakly in $ L^{r}(0,T; L^{\frac{6r}{3r-4}}(\tore)^3) \textrm { for all } r \ge 2, $
while the sequence
 $  q_{\alpha_j} $  converges strongly to $p$ in the space $ L^p([0,T];L^p(\tore)) \textrm{ for all } \frac{4}{3}\le  p<\frac{5}{3}, $ and weakly  in the space  $L^{\frac{r}{2}}(0,T; L^{\frac{3r}{3r-4}}(\tore)) \textrm { for all } r \ge 2.$
\end{Theorem}
\textbf{Proof of Theorem \ref{conmhd}.} As in the proof of Theorem  \ref{1deuxieme} we can show that for all   $ \frac{4}{3} \le  p < \frac{5}{3} $  
we have 
\begin{align}
\overline{\bw_\alpha\otimes\bw_\alpha}  &\rightarrow \bw \otimes \bw &&\textrm{ strongly in  }
L^p(0,T;L^p(\tore)^{3 \times 3}),\\
\overline{\bw_\alpha\otimes\bW_\alpha}  &\rightarrow \bw \otimes \B &&\textrm{ strongly in  }
L^p(0,T;L^p(\tore)^{3 \times 3}),\\
\overline{\bW_\alpha\otimes\bw_\alpha}  &\rightarrow \B \otimes \bw &&\textrm{ strongly in  }
L^p(0,T;L^p(\tore)^{3 \times 3}),\\
\overline{\bW_\alpha\otimes\bW_\alpha}  &\rightarrow \B \otimes \B &&\textrm{ strongly in  }
L^p(0,T;L^p(\tore)^{3 \times 3}).
\end{align}
The above $L^p L^p$ convergences combined with the fact that  $\vec{w}_{\alpha}$ and $\bW_{\alpha} $  belong to  the energy space of the solutions of the Navier-Stokes equations and the  Aubin-Lions compactness  Lemma allow us to take the limit $ \alpha \rightarrow 0$   in order to deduce that $(\bw_{\alpha}, \bW_{\alpha}, q_{\alpha} )$  converge to $ (\bv, \B, p)$  a weak solution to the MHD equations. 
 The rest can be done in exactly
way as in \cite{LL06b}, so we omit the details.\\

 \textbf{Acknowledgement:}
The author thanks professeur R. Lewandowski for interesting discussion
about this paper.


\begin{thebibliography}{10}

\bibitem{adamsstolz}
N.A. Adams and S.~Stolz.
\newblock {\em Deconvolution methods for subgrid-scale approximation in
  large-eddy simulation. In: Modern Simulation Strategies for Turbulent Flow}.
\newblock R.T. Edwards, ed., 2001.

\bibitem{A01}
Hani Ali.
\newblock On a critical {L}eray-$\alpha$ model of turbulence.
\newblock {\em Submitted to DCDS-B}.

\bibitem{HA11}
Hani Ali.
\newblock Mathematical results for some $\displaystyle{\alpha}$ models of
  turbulence with critical and subcritical regularizations.
\newblock {\em Submitted to JMFM}, 2011.

\bibitem{bardina80}
J.~Bardina, J.~Ferziger, and W.~Reynolds.
\newblock Improved subgrid scale models for large eddy simulation.
\newblock {\em American Institute of Aeronautics and Astronautics Paper},
  80:80--1357, 1980.

\bibitem{berselli}
L.C. Berselli, T.~Iliescu, and W.J. Layton.
\newblock {\em Mathematics of Large Eddy Simulation of Turbulent Flows}.
\newblock Springer-Verlag, Berlin, 2006.

\bibitem{bresslilewandowski}
L.C. Berselli and L.~Lewandowski.
\newblock Convergence of approximate deconvolution models to the filtered
  navier-stokes equations.
\newblock {\em under revision in Ann. IHP}, 2011.

\bibitem{CLT06}
Y.~Cao, E.~M. Lunasin, and E.~S. Titi.
\newblock Globall well-posdness of the three dimensional viscous and inviscid
  simplified bardina turbulence models.
\newblock {\em Comm. Math. Sci.}, 4(4):823--848, 2006.

\bibitem{CFHOTW99b}
S.~Chen, C.~Foias, D.~Holm, E.~Olson, E.~S. Titi, and S.~Wynne.
\newblock The {C}amassa-{H}olm equations and turbulence.
\newblock {\em Physica D}, D133:49--65, 1999.

\bibitem{dunca06}
A.~{Dunca} and Y.~{Epshteyn}.
\newblock On the {S}tolz-{A}dams deconvolution model for the large-eddy
  simulation of turbulent flows.
\newblock {\em SIAM J. Math. Anal.}, 37(6):1890--1902, 2006.

\bibitem{DuLi72}
G.~Duvaut and J.~L. Lions.
\newblock {\em Les in\'{e}quations en m\'{e}canique et en physique}.
\newblock Dunod, Paris, 1972.

\bibitem{germano}
M.~Germano.
\newblock Differential filters for the large eddy simulation of turbulent
  flows.
\newblock {\em Phys. Fluids}, 29:1755--1757, 1986.

\bibitem{GH05}
B.~J. Geurts and D.~D. Holm.
\newblock Leray and {L}{A}{N}{S}-alpha modeling of turbulent mixing.
\newblock {\em {J}ournal of {T}urbulence}, 00:1--42, 2005.

\bibitem{HLT10}
M.~Holst, E.~Lunsain, and G.~Tsogtgerel.
\newblock Analysis of a general family of regularized navier-stokes and mhd
  models.
\newblock {\em Journal of Nonlinear Science}, 20(2):523--567, 2010.

\bibitem{latr210}
A.~Labovschii and C.~Trenchea.
\newblock Approximate deconvolution models for magnetohydrodynamics.
\newblock {\em Numerical Functional Analysis and Optimization},
  31(12):1362--1385, 2010.

\bibitem{latr2010}
A.~Labovschii and C.~Trenchea.
\newblock Large eddy simulation for flows.
\newblock {\em Journal of Mathematical Analysis and Applications},
  377(2):516--513, 2011.

\bibitem{LL06}
W.~Layton and R.~Lewandowski.
\newblock Residual stress of approximate deconvolution large eddy simulation
  models of turbulence.
\newblock {\em {J}ournal of {T}urbulence}, 7(46):1--21, 2006.

\bibitem{LL03}
W.~{L}ayton and {R. Lewandowski}.
\newblock A simple and stable scale similarity model for large eddy simulation:
  energy balance and existence of weak solutions.
\newblock {\em Applied Math. letters}, 16:1205--1209, 2003.

\bibitem{LL06b}
W.~{Layton} and {R. Lewandowski}.
\newblock On a well posed turbulence model.
\newblock {\em Continuous Dynamical Systems series B}, 6(1):111--128, 2006.

\bibitem{LT2007}
J.~Linshiz and E.~S. Titi.
\newblock Analytical study of certain magnetohydrodynamic-alpha models.
\newblock {\em J. Math. Phys.}, 48(6), 2007.

\bibitem{RT83mhd}
{M. Sermange} and {R. Temam}.
\newblock {Somes Mathematicals questions related to the MHD equations}.
\newblock {\em Rapport de recherche, INRIA}, (185), 1983.

\bibitem{OT2007}
E.~Olson and E.~S. Titi.
\newblock Viscosity versus vorticity stretching: Global well-posdness for a
  family of navier-stokes-alpha-like models.
\newblock {\em Nonlinear Analysis}, 66:2427--2458, 2007.

\bibitem{sim87}
J.~Simon.
\newblock Compact sets in the spaces $l^p(0,t;b)$.
\newblock {\em Annali di Mat. Pura ed Applic.}, 146:65--96, 1987.

\bibitem{Wa70}
Wolfgang Walter.
\newblock {\em Differential and integral inequalities}.
\newblock Translated from the German by Lisa Rosenblatt and Lawrence Shampine.
  Ergebnisse der Mathematik und ihrer Grenzgebiete, Band 55. Springer-Verlag,
  New York, 1970.

\end{thebibliography}
\end{document}